\documentclass[opre,nonblindrev]{informs3} % current default for manuscript submission
\usepackage{etex}

\usepackage{enumerate}
\DoubleSpacedXI % Made default 4/4/2014 at request
\usepackage{endnotes}
\let\footnote=\endnote
\let\enotesize=\normalsize
\def\notesname{Endnotes}%
\def\makeenmark{$^{\theenmark}$}
\def\enoteformat{\rightskip0pt\leftskip0pt\parindent=1.75em
\leavevmode\llap{\theenmark.\enskip}}

\usepackage{amsmath}
\usepackage{amsfonts}
\usepackage{color}
\usepackage{tikz,pgfplots,pgfplotstable,booktabs}
\usetikzlibrary{shapes,arrows,calc}
\usepackage{caption}
\usepackage{subcaption}
\usepackage{hyperref}
\usepackage{multirow}
\usepackage{empheq}

\usepackage{natbib}
 \bibpunct[, ]{(}{)}{,}{a}{}{,}%
 \def\bibfont{\small}%
\TheoremsNumberedThrough     % Preferred (Theorem 1, Lemma 1, Theorem 2)
\ECRepeatTheorems

\EquationsNumberedThrough    % Default: (1), (2), ...
\newcommand*{\pii}{p_I}
\newcommand*{\pf}{p_F}
\newcommand*{\pw}{p_{W}}
\newcommand*{\pv}{p_V}

\newcommand{\lf}{\lambda^f}
\newcommand{\lb}{\lambda^b(\omega)}
\newcommand{\lbh}{\overline{\lambda^b}}

\newcommand{\pt}{\widetilde{p}}
\newcommand{\ct}{\widetilde{c}}
\newcommand{\zt}{\widetilde{z}}

\newcommand{\ph}{\hat{p}}

\newcommand*{\mi}{\overline{p}_I}
\newcommand*{\mf}{\overline{p}_F}
\newcommand*{\mw}{\overline{p}_W}
\newcommand*{\uw}{\widehat{p}_W}

\newcommand{\ww}{W(\omega)}
\newcommand{\sw}{s(\omega)}
\newcommand{\pfu}{p^{+}_F(\omega)}
\newcommand{\pfd}{p^{-}_F(\omega)}
\newcommand{\iw}{\Delta p_W(\omega)}
\newcommand{\iv}{\Delta p_V(\omega)}

\newcommand*{\ci}{c_I}
\newcommand*{\cf}{c_F}
\newcommand*{\cfu}{c^{+}_F}
\newcommand*{\cfd}{c^{-}_F}

\newcommand{\nvs}{\frac{\cf-\cfd}{\cfu-\cfd}}

\newcommand{\pp}[1]{\left( #1 \right)}

\newcommand*{\om}{\Omega}

\newcommand*{\aom}{\forall \omega \in \Omega}

\newcommand{\ifq}{\text{\quad if \quad}}
\newcommand{\andq}{\text{\quad and \quad}}

\newcommand{\der}[2]{\frac{\partial #1}{\partial #2}}

\newcommand{\inw}[3]{\int_{#1}^{#2} #3 f(\omega)d\omega}
\newcommand{\inww}[3]{\int_{#1}^{#2} #3 F(\omega)d\omega}

\newcommand{\of}[2]{\underset{#1}{{\rm  Minimize}} \quad {#2} }

\newcommand{\acc}[1]{{\rm s.t.} \quad {#1} }
\newcommand{\ac}[1]{\phantom{{\rm s.t.}} \quad {#1} }
\newcommand{\miin}[1]{\min \pp{#1}}
\newcommand{\maax}[1]{\max \pp{#1}}

\newcommand{\f}[1]{f \pp{#1}}
\newcommand{\F}[1]{F \pp{#1}}
\newcommand{\Fi}[1]{F^{-1} \pp{#1}}
\newcommand{\FF}[1]{\pp{1-F\pp{#1}}}



\newcommand{\tred}[1]{\textcolor{red}{#1}}

\begin{document}
%%%%%%%%%%%%%%%%

%\setcounter{page}{1}
%\DOI{10.1287/xxxx.0000.0000}%

\RUNAUTHOR{Morales and Pineda}

\RUNTITLE{On the Inefficiency of the Merit Order in Forward Markets with Uncertain Supply}

% \TITLE{Bundling Information Goods of Decreasing Value}
\TITLE{On the Inefficiency of the Merit Order in Forward Electricity Markets with Uncertain Supply}

\ARTICLEAUTHORS{%
\AUTHOR{Juan Miguel Morales}
%\AFF{Department of Applied Mathematics and Computer Science, Technical University of Denmark, Kgs. Lyngby, Denmark, \EMAIL{jmmgo@dtu.dk; juanmi82mg@gmail.com}} %, \URL{}}
\AFF{Technical University of Denmark, Kgs. Lyngby, Denmark \EMAIL{jmmgo@dtu.dk; juanmi82mg@gmail.com}} %, \URL{}}
\AUTHOR{Salvador Pineda}
%\AFF{Department of Mathematical Sciences, University of Copenhagen, Copenhagen, Denmark, \EMAIL{spinedamorente@gmail.com}}
\AFF{University of Copenhagen, Copenhagen, Denmark, \EMAIL{spinedamorente@gmail.com}}
%\AUTHOR{Marco Zugno}

%\AFF{Department of Applied Mathematics and Computer Science, Technical University of Denmark, Kgs. Lyngby, Denmark, \EMAIL{mazu@dtu.dk}}
%\AFF{Technical University of Denmark, Kgs. Lyngby, Denmark, \EMAIL{mazu@dtu.dk}}
} % end of the block

\ABSTRACT{%
This paper provides insight on the economic inefficiency of the classical merit-order dispatch
in electricity markets with uncertain supply. For this, we
consider a power system whose operation is driven by a two-stage
electricity market, with a forward and a real-time market.  We
analyze two different clearing mechanisms: a
\emph{conventional} one, whereby the forward and the balancing
markets are independently cleared following a merit order, and a
\emph{stochastic} one, whereby both market stages are co-optimized with a view to minimizing the expected aggregate system operating cost.
%to minimize the expectation of the total system
%operating costs.
We first derive analytical formulae to
determine the dispatch rule prompted by the co-optimized two-stage market for a stylized power system with flexible, inflexible and stochastic power generation and infinite transmission capacity. This
exercise %shows that it is not possible for the conventional
%market-clearing mechanism to craft the equivalent  supply
%cost function, unless the market players cooperate and share
%information on their marginal costs. It also
sheds light on the conditions for the stochastic market-clearing
mechanism to break the merit order. We then introduce and
characterize two enhanced variants of the conventional two-stage
market that result in either price-consistent or cost-efficient
merit-order dispatch solutions, respectively. The first of these
variants corresponds to a conventional two-stage market that
allows for virtual bidding, while the second requires that the
stochastic power production be centrally dispatched.  Finally, we
discuss the practical implications of our analytical results and
illustrate our conclusions through examples.
% Enter your abstract
}%
%\KEYWORDS{deterministic inventory theory; infinite linear programming duality;
%  existence of optimal policies; semi-Markov decision process; cyclic schedule}


\maketitle
%%%%%%%%%%%%%%%%%%%%%%%%%%%%%%%%%%%%%%%%%%%%%%%%%%%%%%%%%%%%%%%%%%%%%%

% Samples of sectioning (and labeling) in OPRE
% NOTE: (1) \section and \subsection do NOT end with a period
%       (2) \subsubsection and lower need end punctuation
%       (3) capitalization is as shown (title style).
%
%\section{Introduction.}\label{intro} %%1.
%\subsection{Duality and the Classical EOQ Problem.}\label{class-EOQ} %% 1.1.
%\subsection{Outline.}\label{outline1} %% 1.2.
%\subsubsection{Cyclic Schedules for the General Deterministic SMDP.}
%  \label{cyclic-schedules} %% 1.2.1
%\section{Problem Description.}\label{problemdescription} %% 2.

% Text of your paper here

\section{Introduction}
\label{Intro}

% Exchange/pool definition, dispatch, pricing
Electricity markets are typically arranged as sequences of exchanges or pools, where producers and possibly consumers submit offers and bids specifying the amount of electricity they are willing to deliver to or withdraw from the network and at what unit price. During the market-clearing process, a market operator determines the optimal dispatch by accepting a subset of the submitted offers for electricity production and, possibly, bids for consumption. Furthermore, a market price or a set of prices, one for each node across the electricity network, is set so that the dispatched generation and consumption blocks are profitable according to the respective offer and bid prices.

% Why different pools
Generally, electricity markets comprise different floors for trading electricity arranged in a sequential fashion up to the time when electricity is delivered. The presence of a real-time or balancing market is necessary as electricity is a non-storable commodity, and since imbalances between supply and offtake result in deviations of the system frequency that can harm machines and appliances connected to the network. On the other hand, market floors that clear hours ahead of electricity delivery are also needed to guarantee the participation of units with slower response time, e.g., nuclear power plants. Although market structures and regulations vary significantly from country to country, a common trait is the presence of at least a day-ahead market stage (besides a real-time one) clearing around noon on the day prior to the delivery of electricity. In many European countries, day-ahead markets account for the bulk of the total trading of electricity \citep{Weber2010}. Although other market stages may exist, we consider in the following a two-stage market comprising a forward (day-ahead) and a real-time floor as an abstraction of electricity markets comprising different market stages.

% Sequential market clearing: merit-order definition
Traditionally, electricity markets are cleared in a sequential manner and independently from one another, i.e., without considering the impact of the market-clearing decision on the operation of future markets stages. In such a \emph{conventional} arrangement, each market stage is cleared by dispatching producers and consumers in order according to increasing marginal costs and decreasing marginal utilities, respectively. Under the assumption that producers and consumers bid their ``true'' marginal costs and benefits, such a dispatch of the forward market based on the \emph{merit-order} principle results in schedules that maximize the social welfare \emph{exclusively} for that market stage. However, there is no guarantee that it maximizes the total social welfare, i.e., aggregated across sequential market stages, in the long run. Indeed, the application of the merit-order principle may leave fast-ramping units out of the forward dispatch and potentially result in a lack of flexibility, and hence inefficiency, at the real-time stage.

% Impact of renewables
Arguably, the suboptimality of a merit-order based forward dispatch is exacerbated by the increasing penetration of partly-predictable renewables such as wind and solar. These generation sources are characterized by a zero (or near-zero) production cost per unit and their dispatch is often prioritized by market operators. As a result, owners of renewable production facilities submit price-inelastic forward offers, where a certain quantity, typically a point forecast of future generation, is offered at zero price \citep{Morales2014integrating}. In a merit-order based forward dispatch, a large penetration of renewables may, besides pushing flexible resources out of the market, increase the need for flexibility at the balancing stage, which results from their uncertain nature.

% Beyond sequential market clearing: SUC, stochastic dispatch
The attention received by the impact of uncertain renewable generation on electricity markets has grown in the technical literature along with their actual deployment in power systems. Numerous papers, see for example \cite{Jonsson2010}, have assessed the downward pressure exercised by uncertain renewable power on electricity market prices, i.e., the so-called \emph{merit-order effect}. Among the most notable consequences of the merit-order effect is the emergence of negative prices, which are caused by the combination of large volumes of zero-price offers from renewable suppliers and thermal producers that are willing to incur occasional losses to avoid wear-and-tear of their units. \cite{Hildmann2014} claims that the removal of feed-in-tariff schemes would partially solve this problem by inciting owners of renewable generation facilities to internalize forecasting-error (balancing) costs in their offer, hence resulting in positive marginal costs of renewables. However, balancing costs are hard to estimate and vary from hour to hour.

Other works have focused on centralized solutions for the
coordinated clearing of the forward and real-time market stages
\citep{Pritchard2010,morales2012pricing}. In these proposals, the
forward dispatch is determined by a market operator that makes use
of two-stage stochastic programming to account for the balancing
costs due to the uncertain supply and minimize the expected
aggregate system cost. While stochastic dispatch models do not
comply with the merit-order at the forward stage, they guarantee
both revenue adequacy for the market operator
\citep{Pritchard2010} and cost-recovery for the producers
\citep{morales2012pricing} in expectation. Furthermore, they
improve the long-run social welfare. An alternative solution where
the merit-order at the forward market is only broken for providers
of uncertain supply is proposed in \cite{morales2014electricity}.
That work shows that part of the improvement in expected social
welfare achievable by the stochastic dispatch can be captured by a
classical forward market where the dispatch of uncertain renewable
generators is forced by the market operator. In general, though,
the problem of determining the optimal forward offer or dispatch
for stochastic power suppliers has no closed form solution in
terms of a simple statistic of the forecast distribution of
uncertain supply. \citet{Zavala2015} investigates stochastic
dispatch schemes with $\ell_1$ and $\ell_2$ penalization of
imbalances where the optimal dispatch converges to the conditional
median and mean of the uncertain supply. However, a generalization
of these properties is not possible, as the underlying penalty
assumptions may clash with the offering preferences of the
flexible producers, which ultimately drive the penalization of the
imbalances at the real-time market.

This paper aims to shed light on the impact that preserving the
merit order in forward electricity markets with uncertain supply
may have on market efficiency. By \emph{market efficiency} we mean
the ability of the market to minimize the expected aggregate
system operating cost (cost efficiency) \emph{and} to deliver a
set of prices such that the forward price equals the expectation
of the real-time price (price consistency) --- a number of authors
have underlined the benefits of price-consistent market settings
\citep{bessembinder2002equilibrium, kaye1990forward, Zavala2015}.
For this purpose, we construct mathematical models for four
different types of two-stage markets, namely:
\begin{enumerate}
    \item A market that is price-consistent \emph{and} guarantees
    maximum cost-efficiency. A market with such properties is
    achieved by co-optimizing the forward dispatch and the real-time re-dispatch
    through the use of stochastic programming. We show that this
    market produces, however, dispatch solutions that break the
    merit order. Furthermore, the practical implementation of this
    market is not without its challenges and problems in terms of
    revenue adequacy, cost recovery, arbitrariness in the
    probabilistic characterization of the uncertain supply, etc.
    \citep{morales2014electricity}. Therefore, we just use it here as an ``ideal''
    benchmark that theoretically achieves the highest market
    efficiency.
    \item A market that follows the merit order, but that is, in
    general, price inconsistent and cost-inefficient. This is the
    case of a conventional two-stage market, where the forward and
    the real-time settlements are \emph{not} co-optimized, and
    where the uncertain supply is systematically dispatched to a certain statistic of its forecast probability distribution (most commonly, the conditional
    expectation).
    \item A price-consistent market that preserves the merit
    order. We construct this market from the conventional two-stage market described in point 2 above,
    by introducing a risk-neutral virtual bidder that arbitrages between the forward and the real-time markets. We show that, in order to ensure price consistency, such a market
    may have to give up on cost efficiency.
    \item A market that renders the most cost-efficient dispatch
    among those that respect the merit order. The practical
    translation of this market is that of a conventional two-stage
    market, like the one described in point 2 above, in which the uncertain supply is centrally dispatched by a non-profit, all-knowing organization such as
    an Independent System Operator. We show that, in order for
    this market to ensure maximum cost-efficiency, while complying
    with the merit order, it may have to give up on price
    consistency.
 \end{enumerate}

Unlike other works that rely on computational simulation for their analysis, e.g.,
\citet{Bouffard08, khazaei2014effects, morales2012pricing, morales2014electricity}, we derive closed-form solutions to the mathematical models
describing these four types of markets for a stylized power system with flexible, inflexible and stochastic power generation and infinite transmission capacity. This exercise allows us to
characterize the dispatch solutions prompted by these markets and
identify conditions for their equivalence or dissimilarity. We
accompany this analytical insight with a meaningful discussion on
the practical implications of our results and illustrate our main
conclusions through examples.

The structure of the paper is the following. Sections
\ref{Conventional} and \ref{Stochastic} deal with markets 1 and 2,
respectively. More specifically, we define and formulate the
conventional and the stochastic two-stage electricity markets for
a stylized power system and provide the closed forms of the
dispatch rules that each of these markets induce. In these two
sections we also introduce some important concepts that are
repeatedly used throughout the paper. In Section
\ref{VirtualBidding}, we focus on market 3 and provide conditions
under which the conventional and the stochastic market-clearing
models are equivalent in the case that virtual bidding is allowed. Section
\ref{Centralized} deals with market 4, that is, with the case of a
conventional two-stage market in which the stochastic power
production is centrally dispatched with the aim of minimizing the
expected aggregate system operating cost. The study of this market
setting allows us to identify conditions under which a
conventional two-stage market, even if price-consistent, does not
deliver the most cost-efficient merit-order dispatch. Finally,
conclusions are drawn in Section \ref{Conclusion}.

%We consider a power system whose operation is driven by a two-stage electricity market, which consists of a forward (typically, day-ahead) and a real-time or balancing market. The forward market, also known as \emph{spot market}, is normally cleared one day before the physical delivery of the traded electricity is to take place. Electricity transactions in the balancing market are, on the contrary, for their (almost) immediate delivery.

%The forward market provides slow generating units---i.e., the \emph{inflexible} power plants---with enough time to reliably set their power production profile based on the market outcomes. On the other hand, the balancing market allows stochastic power units, whose power production can only be predicted with limited accuracy at the time the forward market is cleared, to trade their energy imbalances with the more \emph{flexible} power plants.

%Next we briefly introduce two different mechanisms for clearing the two-stage market, which essentially differ in whether the operation of the forward and the balancing market is co-optimized or not. This, in turn, will determine how proficient the two-stage market is in coping with the power imbalances brought into the system by the stochastic power producers.

\section{Conventional or Inefficient Two-stage Market (ConvM)} \label{Conventional}

%The stochastic and the conventional two-stage markets essentially differ in whether the operation of the forward and the balancing markets is co-optimized or not, respectively. This, in turn, will determine how proficient the two-stage market is in coping with the power imbalances brought into the system by the stochastic power producers. \tred{Remove this paragraph?}

We first formulate the model of the conventional two-stage market, where the operation of the forward and the balancing markets is \emph{not} co-optimized. Each market, therefore, attempts to minimize operating costs independently. We describe such a two-stage market as \emph{inefficient}, because it results in higher operating costs in the long run.

The market models that we introduce throughout the paper are all tailored to the stylized power system described below.

\begin{definition}[Stylized power system]\label{PowerSystem}
Our stylized power system has infinite transmission capacity and consists of inflexible, flexible and stochastic power generation technologies with capacities $\mi > 0$, $\mf > 0$, and $\mw > 0$, in that order. We denote the marginal cost of the flexible and inflexible generating capacity by $\cf > 0$ and $\ci > 0$, respectively, and assume that the marginal cost of the stochastic power production is zero. Furthermore, we use $\cfu$ to represent the incremental cost incurred by the flexible generating capacity for marginally increasing its production for balancing (upward regulation) and $\cfd$ to denote the incremental utility that it obtains from marginally decreasing its production for
balancing (downward regulation). The demand $l$ is inelastic and known with certainty with a cost of involuntary curtailment denoted by $v$. In this stylized power system, it holds that $v > \cfu \geq \cf \geq \cfd \geq 0$. Finally, the power production from renewable sources is characterized as a random variable $W$ defined on some probability space $\left(\Omega, \mathcal{F}, P\right)$. Let $F(\cdot)$ and $f(\cdot)$ denote the cumulative distribution function and the probability density function of $W$, respectively. To comfortably deal with the case $F(0) > 0$ in our mathematical derivations, we use the generalized inverse cumulative distribution function $F^{-1}(\alpha) = \text{inf}\{x \in \mathbb{R}: F(x) \geq \alpha\}$.
\end{definition}

The aim of the forward market is to determine the dispatch of flexible power producers ($\pf$), inflexible power producers ($\pii$), and stochastic power producers ($\pw$) that minimizes (forward) system operating costs, that is,
\begin{subequations}\label{Ineff}
\begin{align}
& \underset{\pf, \pii, \pw}{{\rm Minimize}}\quad \ci \pii + \cf \pf\\
& {\rm s.t.}\quad \pf + \pii + \pw - l = 0: \lambda^{f}\;, \label{Ineff:PB}\\
& \phantom{s.t.}\quad 0 \leq \pf \leq \mf\;, \label{Ineff:LimitsF}\\
& \phantom{s.t.}\quad 0 \leq \pii \leq \mi\;, \label{Ineff:LimitsI}\\
& \phantom{s.t.}\quad 0 \leq \pw \leq \uw \label{Ineff:LimitsW}
%
%& \phantom{s.t.}\quad \pf, \pii  \geq 0 \label{Ineff:Pos}
%
\end{align}
\end{subequations}
%
%where $\delta^{f}$ is the vector of state variables that determine power flows in the meshed power system. For example, $\delta^{f}$ may represent the vector of voltage angles if the algorithm for clearing the forward market includes a DC power flow model of the transmission network.

Equation~\eqref{Ineff:PB} enforces the power balance.
%For ease of exposition, and without loss of generality, we assume that demand $l$ is inelastic and known with certainty (\tred{remove since included in Definition 1}).
We denote the Lagrange multiplier associated with the power balance equation by $\lambda^{f}$, which defines the forward electricity price. %, which are used to price electricity in many markets.
The set of inequalities~\eqref{Ineff:LimitsF}--\eqref{Ineff:LimitsW} imposes upper and lower bounds on the dispatch of the different power producers. We explicitly indicate parameter $\uw$, which stands for the power production that is expected from the stochastic power producers, in~\eqref{Ineff:LimitsW}. We do so to note that, typically, the amount of stochastic power production that can be cleared in the forward market is capped to this expectation \citep{Bouffard08, Cadre2014,Oggioni2014,Zavala2015}.

The following proposition provides the optimal solution to problem~\eqref{Ineff}.

\begin{proposition}[The merit-order dispatch solution]\label{Prop:MeritOrder}
Consider the stylized power system described in Definition~\ref{PowerSystem}, where, in addition, it holds that $\ci < \cf$. Optimization problem~\eqref{Ineff} prompts the following dispatch rule:

\vspace{4mm}
\begin{center}
{\footnotesize
\begin{tabular}{|l|lll|l|}
\hline
Rule \# & $\pw$ & $\pii$ & $\pf$ & applies if \\
\hline
1. & $l$ & $0$ & $0$ &  $0 \leq l \leq \uw$ \\
2. & $\uw$ & $l-\uw$ & $0$ &  $\uw < l \leq \uw+\mi$ \\
3. & $\uw$ & $\mi$ & $l-\uw-\mi$ & $\uw+\mi < l \leq \uw+\mi+\mf$ \\
4. & \multicolumn{3}{c|}{infeasible} & $\uw+\mi+\mf < l$ \\
\hline
\end{tabular}}\end{center}
\end{proposition}

This result is well known (see, e.g., \citet[Chapter 5]{gomez2008electric}) and therefore, its proof is omitted here.

Note that values of load $l > \uw+\mi+\mf$ render problem~\eqref{Ineff} infeasible, because we have not considered the possibility of shedding load in the forward market.

The conventional market produces a forward dispatch whereby power production is cleared following the so-called \emph{merit order}, i.e., the power plants with the lowest marginal costs are dispatched first, until the system demand is satisfied. In this case, the dispatch solution prompted by this market follows from the intersection of the (forward) supply cost function of the system with the marginal utility demand curve. The supply cost function is built by sorting the marginal cost functions of the individual power plants in increasing order. Consequently, all the production with a marginal cost lower than the one determined by the said intersection is dispatched.

Optimization problem~\eqref{Ineff} results, therefore, in the forward dispatch $\left(\pf, \pii, \pw\right)$ given by Proposition~\ref{Prop:MeritOrder}. Since the stochastic power production cannot be perfectly predicted, the (random) power imbalance $W-\pw$ is to be covered in the balancing market. For this purpose, flexible power plants can be re-dispatched and/or the amount of load and stochastic power production can be curtailed.

Suppose a specific realization $W(\omega)$, $\omega \in \Omega$, of the stochastic power production. The balancing market determines the most economical vector of re-dispatch actions that accommodates the power imbalance $W(\omega)-p_{W}$, that is,
%
%\begin{subequations}\label{IneffBM}
%%
%\begin{align}
%%
%& \underset{r_{\omega}, \delta^{b}_{\omega}}{{\rm Minimize}}\quad \mathcal{C}^{b}\left(r_{\omega}\right)\\
%%
%& {\rm s.t.}\quad h^{b}\left(r_{\omega}, \delta^{b}_{\omega}, \delta^{f \star}\right) + W_{\omega}-p_{W}^{\star} = 0: \lambda^{b}_{\omega}\;, \label{IneffBM:PB}\\
%%
%& \phantom{s.t.}\quad g^{b}\left(r_{\omega}, \delta^{b}_{\omega}, p_{F}^{\star}, p_{W}^{\star}; W_{\omega}\right) \leq 0\;, \label{IneffBM:Limits}
%%
%\end{align}
%%
%\end{subequations}
%
\begin{subequations}\label{IneffBM}
\begin{align}
& \of{\pfu, \pfd, \iw, \sw}{v \sw + \cfu \pfu - \cfd \pfd} \\
& \acc{\sw + \pfu - \pfd + \iw = 0: \lambda^{b}(\omega)\;, } \label{IneffBM:PB}\\
& \ac{0 \leq \pfd \leq \pf\;, } \label{IneffBM:LimitsDW}\\
& \ac{0 \leq \pfu \leq \mf - \pf\;, \label{IneffBM:LimitsUP}} \\
& \ac{0 \leq \pw + \iw \leq \ww\;, \label{IneffBM:LimitsWS}} \\
& \ac{0 \leq \sw \leq l \;, \label{IneffBM:LimitsLS}}
%& \ac{\pfu, \pfd, \sw \geq 0\;.}
\end{align}
\end{subequations}

%
%where $\mathcal{C}^{b}\left(\cdot\right)$ is a real-valued function providing the cost of the balancing actions $r_{\omega}$ and $\delta^{b}_{\omega}$ is the vector of state variables determining the power flows after the re-dispatch of the system.
Equation~\eqref{IneffBM:PB} ensures that the power system is brought to balance by deploying upward or downward regulation from the flexible power unit, i.e., $\pfu$ or $\pfd$, respectively; curtailing load $\sw$ and/or curtailing stochastic power production, which is given by $\ww - \pw - \iw$. The Lagrange multiplier $\lambda^{b}(\omega)$ defines the marginal price that clears the balancing market. The family of inequalities~\eqref{IneffBM:LimitsDW}--\eqref{IneffBM:LimitsLS} set limits on the amount of downward and upward regulation that the flexible power unit can provide, \eqref{IneffBM:LimitsDW} and \eqref{IneffBM:LimitsUP}, respectively; the amount of stochastic power production that can be curtailed \eqref{IneffBM:LimitsWS}, and the amount of load that can be shed \eqref{IneffBM:LimitsLS}.

%\tred{I have removed this sentence: The total cost of operating the power system under the realization $W(\omega)$ of the stochastic power production is, thus, given by $\ci \pii + \cf \pf + v \sw + \cfu \pfu - \cfd \pfd$.}

The total expected cost of operating the power system is given by $\ci \pii + \cf \pf + \mathcal{C}^b(\pii,\pf,\pw)$, where $\mathcal{C}^b(\pii,\pf,\pw)$ represents the expected balancing cost computed in the proposition below.

\begin{proposition}[Expected balancing cost]\label{Prop:Cost2Stage}
Consider the stylized power system described in Definition~\ref{PowerSystem} with given forward dispatch quantities $\pii,\pf,\pw$. The expected balancing cost $\mathcal{C}^b(\pii,\pf,\pw)$ is computed as
\begin{equation}
\mathcal{C}^b(\pii,\pf,\pw) = v\inww{0}{\pf+\pw-\mf}{} + \cfu\inww{\pf+\pw-\mf}{\pw}{} + \cfd\inww{\pw}{\pf+\pw}{} - \cfd\pf
\end{equation}
\end{proposition}

The proof of this proposition is included in Appendix \ref{Proof:Cost2Stage}. Note that the expected balancing cost is actually independent of $\pii$.

\section{Stochastic or Efficient Two-stage Market (StoM)} \label{Stochastic}
We now build the model of a two-stage market where the operation of the forward and the balancing markets is co-optimized. To this aim, one just needs to replace optimization problem~\eqref{Ineff}
with an alternative market-clearing mechanism that seeks to minimize the \emph{expected} total system operating cost, namely:
%
%\begin{subequations}\label{Eff}
%%
%\begin{align}
%%
%& \underset{p_{F}, p_{I}, p_{W},\delta^{f}; r(\omega), \delta^{b}(\omega)}{{\rm Minimize}}\quad \mathcal{C}^{f}\left(p_{F}, p_{I}, p_{W}\right) + E_{\Omega}\left[\mathcal{C}^{b}\left(r(\omega)\right)\right]\\
%%
%& {\rm s.t.}\quad h^{f}\left(p_{F}, p_{I}, p_{W},\delta^{f}\right) - l = 0: \nu^{f}\;, \label{Eff:PB}\\
%%
%& \phantom{s.t.}\quad g^{f}\left(p_{F}, p_{I}, p_{W}, \delta^{f}\right) \leq 0\;, \label{Eff:Limits}\\
%%
%& \phantom{s.t.}\quad h^{b}\left(r(\omega), \delta^{b}(\omega), \delta^{f}\right) + W(\omega)-p_{W} = 0\;, \forall \omega \in \Omega \label{EffBM:PB}\\
%%
%& \phantom{s.t.}\quad g^{b}\left(r(\omega), \delta^{b}(\omega), p_{F}, p_{W}; W(\omega)\right) \leq 0\;, \forall \omega \in \Omega \label{EffBM:Limits}
%%
%\end{align}
%%
%\end{subequations}
%%
\begin{subequations}\label{Eff}
\begin{align}
& \of{\pii, \pf, \pw; \pfu, \pfd, \iw, \sw}{\ci \pii + \cf \pf + \inw{\om}{}{\pp{v \sw + \cfu \pfu - \cfd \pfd}}} \\
& \acc{\pf + \pii + \pw - l = 0: \nu^{f}}\\
& \ac{\sw + \pfu - \pfd + \iw = 0: \nu^{b}(\omega) f(\omega), \quad \aom} \label{EffBM:PB}\\
& \ac{0 \leq \pf \leq \mf}\\
& \ac{0 \leq \pii \leq \mi}\\
& \ac{0 \leq \pfd \leq \pf, \quad \aom} \label{EffBM:LimitsDW}\\
& \ac{0 \leq \pfu \leq \mf - \pf, \quad \aom} \label{EffBM:LimitsUP}\\
& \ac{0 \leq \pw + \iw \leq \ww: \underline{\gamma}(\omega), \overline{\gamma}(\omega), \quad \aom} \label{EffBM:LimitsWS}\\
& \ac{0 \leq \sw \leq l, \quad \aom}\label{EffBM:LimitsLS}
%& \ac{\pii, \pf \geq 0}\\
%& \ac{\pfu, \pfd, \sw \geq 0, \quad \aom}
\end{align}
\end{subequations}

where the expectation of the balancing cost $\inw{\om}{}{\pp{v \sw + \cfu \pfu - \cfd \pfd}}$ is taken over the probability space $\left(\Omega, \mathcal{F}, P\right)$ on which the stochastic power production $W$ is defined. Note that, for ease of notation, we write $y(\omega)$ instead of $y\left(W(\omega)\right)$.

Problem~\eqref{Eff} computes the optimal forward dispatch $\left(\pf, \pii, \pw\right)$ by taking into account the potential cost of the subsequent re-dispatch of the system that is induced by the random power imbalance $W-\pw$. Ideally, problem~\eqref{Eff} also provides the optimal re-dispatch rule $\left(\pfu, \pfd, \iw, \sw\right)$ that guarantees, by enforcing \eqref{EffBM:PB} and \eqref{EffBM:LimitsDW}--\eqref{EffBM:LimitsLS}, the power balance for any possible outcome $W(\omega)$ of the random variable W.

%However, problem~\eqref{Eff} is an infinite programming problem, which contains an infinite number of constraints of the type of \eqref{EffBM:PB} and \eqref{EffBM:LimitsDW}--\eqref{EffBM:LimitsLS} and optimizes functions $\pfu, \pfd, \iw$, and $\sw$. To make this problem tractable, it is customary to assume that $\Omega$ is a countable sample space, formed by a finite number of atoms, outcomes, or scenarios $W(\omega), \omega = 1, \ldots, N_{\Omega}$, each with a probability $\pi(\omega)$ such that $\sum_{\omega = 1}^{N_{\Omega}}{\pi(\omega)} = 1$. \tred{This may confuse the reader since in this paper we are actually solving that problem without relying on scenarios...}

We note that, unlike \eqref{Ineff}, the forward market \eqref{Eff} does not need to arbitrarily cap the dispatch $\pw$ of stochastic power production, because in  \eqref{Eff} the optimization of $\pw$ is driven by the probabilistic characterization of random variable $W$, which is naturally bounded and nonnegative. By the same token, we do not need to impose that $\pw \geq 0$, since $\cfu \geq \cf \geq \cfd \geq 0$.

Finally, the balancing stage of this market is also modeled by \eqref{IneffBM}, but with the optimal forward dispatch $\left(\pf, \pii, \pw\right)$ given by \eqref{Eff}.

We now define some relevant concepts that will be used in the remaining part of this paper.

\begin{definition}\label{Def:price_consistency}
A two-stage market such as~\eqref{Ineff}--\eqref{IneffBM} and~\eqref{Eff}--\eqref{IneffBM} is said to be \emph{price consistent} if the forward price is equal to the expected value of the real-time price, that is,
\begin{equation}
  \lambda^f= E_{\Omega}\left[\lambda^{b}(\omega)\right] = \int_{\Omega}{\lambda^{b}(\omega) f(\omega) d\omega}
\end{equation}
in the case of the conventional two-stage market~\eqref{Ineff}--\eqref{IneffBM}, and
\begin{equation}
  \nu^f= E_{\Omega}\left[\nu^{b}(\omega)\right] = \int_{\Omega}{\nu^{b}(\omega) f(\omega) d\omega}
\end{equation}
in the case of the co-optimized two-stage market~\eqref{Eff}--\eqref{IneffBM}.
\end{definition}

We have taken the term \emph{price consistency} from~\cite{Zavala2015}. Definition~\ref{Def:price_consistency} allows us to formulate the following proposition.

\begin{proposition}[Price consistency of the stochastic two-stage market]\label{Prop:PriceCons}
 The two-stage stochastic market~\eqref{Eff}--\eqref{IneffBM} is price consistent.
\end{proposition}

\proof{Proof.} The optimality conditions of problem \eqref{Eff} imply that the derivative of the Lagrangian with respect to $\pw$ and $\iw$ are equal to 0 at the optimal solution, i.e.,
\begin{equation}
\begin{rcases}
\der{\mathcal{L}}{\pw} = 0 \implies \nu^f + \int_\Omega \left( -\underline{\gamma}(\omega) +  \overline{\gamma}(\omega) \right)d\omega = 0 \\
\der{\mathcal{L}}{\iw} = 0 \implies \nu^b(\omega)f(\omega) - \underline{\gamma}(\omega) +  \overline{\gamma}(\omega) = 0, \forall \omega \in \Omega \\
\end{rcases} \implies \nu^f= \int_{\Omega}{\nu^{b}(\omega) f(\omega) d\omega}
\end{equation}

where $\underline{\gamma}(\omega)$ and $\overline{\gamma}(\omega)$ are the dual variables of constraints \eqref{EffBM:LimitsWS}. \hfill \Halmos \endproof

%\begin{definition}
%Let $z(l)$ be the expected aggregate system operating cost at system load level~$l$, that is, $z(l) = \ci \pii + \cf \pf + \inw{\om}{}{\pp{v \sw + \cfu \pfu - \cfd \pfd}}$ where $\left(p_{F}, p_{I}, \pfu, \pfd, \sw\right)$ is solution to~\eqref{Ineff}--\eqref{IneffBM} (or alternatively, solution to~\eqref{Eff}--\eqref{IneffBM}). The \emph{marginal expected system cost} at system load level $l$, denoted by $\frac{dz}{dl}(l)$, is the incremental expected cost incurred by the power system as a result of a marginal change in the forward system load $l$ under the two-stage market~\eqref{Ineff}--\eqref{IneffBM} (or alternatively, under the two-stage market~\eqref{Eff}--\eqref{IneffBM}). \tred{Do we need this definition?}
%\end{definition}

\begin{definition}
Consider the stochastic market-clearing mechanism~\eqref{Eff}. This mechanism is said to \emph{break} or \emph{violate} the merit order when a conventional generating unit is dispatched in the forward market ahead of a power generating unit with a lower marginal cost.
\end{definition}
Note that this latter definition only applies among conventional power plants and leaves out stochastic power generating units. In the particular case of the stylized power system introduced in Definition~\ref{PowerSystem}, the merit order is violated if $\pii < \mi$ and $\pf > 0$ (with $\ci < \cf$) or if $\pf < \mf$ and $\pii > 0$ (with $\cf < \ci$).

\citet{morales2012pricing} and \citet{morales2014electricity} provide examples of situations in which the generation dispatch $\left(p_{F}, p_{I}, p_{W}\right)$ that is solution to the stochastic market-clearing mechanism~\eqref{Eff} violates the merit order. In this paper, we aim at providing analytical insight on the conditions under which this occurs. Our analysis relies on the following theorem, which is indeed one of the major results of our work.

\begin{theorem}[The stochastic dispatch rule]\label{Prop:ESCF}
Consider the stylized power system described in Definition~\ref{PowerSystem}.
%
%Then, the equivalent supply cost function of this power
%system, $ESCF(x)$, is given by:
%%
%\begin{equation}\label{ESCF}
%ESCF(x) = min\left(c_{I}, \enskip c_{F}^{+}F(x), \enskip c_{F} -
%c_{F}^{-}\left(1-F(x)\right)\right)\;,
%\end{equation}
%where $F(x) = P\left(W \leq x\right)$ is the cdf that characterizes the stochastic power
%production W.
%
Now define the constants
\begin{align}
l_1 &:= F^{-1}\left(\frac{c_F-c_{F}^-}{c_F^{+}-c_F^{-}}\right); \label{L1}\\
l_2 &:= \underset{l \geq 0}{\rm min} \enskip l: \left(v-c_F^+\right)F\left(l-\mf\right) + c_F^+F(l) \geq c_I \label{L2}\\
l_3 &:= \underset{l \geq l_1}{\rm min} \enskip l: c_F + \left(v-c_F^+\right)F\left(l-\mf\right) - c_F^- \left(1-F(l)\right) \geq c_I; \label{L3} \\
l_4 &:= \underset{l \geq l_1 + \mf}{\rm min} \enskip l: \left(v-c_F^-\right)F\left(l-\mf\right) + c_F^- F(l) \geq c_I
\end{align}
to which we assign an infinite value in those cases where the corresponding minimization problem is infeasible. %Note that if $\cfu > \ci$, then $l_2 < \mw$.
The efficient two-stage market~\eqref{Eff} prompts the following dispatch rule, which we refer to as the \emph{stochastic dispatch rule}:

\vspace{4mm} \begin{center}
{\footnotesize
\begin{tabular}{|l|lll|l|l|l|l|}
  \hline
  Rule \# & $p_{W}$ & $p_{I}$ & $p_{F}$ & \multicolumn{4}{c|}{applies if}   \\
  \hline
 1. & $l$ & $0$ & $0$                                                             & \multirow{5}{1.2cm}{$l_1 \geq l_2$}  & \multicolumn{3}{l|}{$0 \leq l \leq l_2$} \\
 2. & $l_2$ &$l - l_2$ & $0$                                                       &                                      & \multicolumn{3}{l|}{$l_2 < l \leq \mi + l_2$}\\
 3. & $l-\mi$ & $\mi$ & $0$                                 &                                      & \multicolumn{3}{l|}{$\mi + l_2 < l \leq \mi + l_1$}\\
 4. & $l_1$ & $\mi$ &$l - l_1 -\mi$                        &                                      & \multicolumn{3}{l|}{$\mi + l_1 < l \leq \mf + \mi + l_1$}\\
 5. & $l - \mf - \mi$ &$\mi$& $\mf$ &                                      & \multicolumn{3}{l|}{$l > \mf + \mi + l_1$}\\
\hline
 6. & $l$ & $0$ & $0$                                                             & \multirow{9}{1.2cm}{$l_1 < l_2$}     & \multicolumn{3}{l|}{$0 \leq l \leq l_1$}                                  \\
 \cline{6-8}
 7. & $l_1$ &$0$ & $l-l_1$                                                       &                                      & \multirow{4}{3cm}{$l_3 \leq l_1 + \mf$} & \multicolumn{2}{l|}{$l_1 \leq l \leq l_3$} \\
 8. & $l_1$ &$l - l_3$ & $l_3 - l_1$                                             &                                      &                                                      & \multicolumn{2}{l|}{$l_3 < l \leq \mi + l_3$} \\
 9. & $l_1$ & $\mi$ & $l - l_1 - \mi$                       &                                      &                                                      & \multicolumn{2}{l|}{$l_3 + \mi < l \leq \mi + \mf + l_1$} \\
10. & $l - \mi -\mf$ & $\mi$ & $\mf$  &                                      &                                                      & \multicolumn{2}{l|}{$l > \mi + \mf + l_1$} \\
\cline{1-4}\cline{6-8}
11. & $l_1$ & $0$ & $l-l_1$                                                       &                                      & \multirow{5}{3cm}{$l_3 > l_1 + \mf$}      & \multicolumn{2}{l|}{$l_1 \leq l \leq \mf + l_1$} \\
\cline{7-8}
12. & $l - \mf$ & $0$ & $\mf$                               &                                      &                                                      & \multirow{3}{1.2cm}{$\exists l_4$ $(c_I \leq v)$}   &  $\mf + l_1 \leq l \leq l_4$\\
13. & $l_4 - \mf$ & $l - l_4$ & $\mf$                       &                                      &                                                      &                       & $l_4 \leq l \leq l_4 + \mi$\\
14. & $l - \mf - \mi$ & $\mi$ & $\mf$ &                                      &                                                      &                      & $l > l_4 + \mi$\\
\cline{7-8}
15. & $l - \mf$ & $0$ & $\mf$                               &                                      &                                                      &  \multirow{1}{1.2cm}{$\nexists l_4$}      &  $l > \mf + l_1$\\
  \hline
\end{tabular}} \end{center}
\vspace{4mm}
\end{theorem}

The proof of this theorem is given in Appendix~\ref{Proof:ESCF}.

%

%, where we also provide closed form expressions for the total and the marginal expected system cost, i.e., $z(l)$ and $\frac{dz}{dl}(l)$, respectively. \tred{This is not included in the Appendix. So I think we should remove this sentence.}

Theorem~\ref{Prop:ESCF} can be used to identify cases in which the stochastic market-clearing mechanism~\eqref{Eff} breaks the merit order. For example, consider a power system whose ability to cope with uncertain supply stems from a small amount of flexible power capacity that is comparatively cheap ($c_{F} < c_{I} \ll v$) and that can provide upward regulation at almost no extra cost ($c_{F} \lesssim c_{F}^{+}$). Under these conditions, the block of dispatch rules 1--5 ($l_1 \geq l_2$) in Theorem~\ref{Prop:ESCF} may apply, resulting in the inflexible generation technology being dispatched over the cheaper flexible one. Intuitively, it makes economic sense, under these circumstances, to withhold flexible capacity from the forward market to have it available for upward regulation in real time at almost no extra cost (bear in mind that the alternative would be to curtail load at the very high cost $v$).

The mirrored case would be that of a power system with a large amount of flexible, but comparatively expensive power capacity ($c_{I} < c_{F} \ll v$) that is able to provide downward regulation at nearly no extra cost ($c_{F}^{-} \lesssim c_{F}$). In this situation, the block of dispatch rules 6--15 ($l_1 < l_2$) would hold and the flexible capacity would be dispatched over the cheaper inflexible one. The intuition here is that it is profitable to commit flexible capacity in the forward market to have it available for accommodating surplus of stochastic power production in real time at nearly no extra cost.

In the sequel, we focus our analysis to those cases in which $c_{I} < c_{F}$. This is a fairly common characteristic of many power systems, where, for instance, the base load is mainly supplied by nuclear and large coal-fired power plants, while balancing is mostly provided by gas-fired power units. Notice that the fact that $c_{I} < c_{F}$ precludes the stochastic dispatch rule 15 in Theorem~\ref{Prop:ESCF}. Next we provide results that are directly derived from Theorem~\ref{Prop:ESCF} for relevant special cases. We state these results in the form of corollaries.

\begin{corollary}[Free downward regulation]\label{Corr:FreeDW}
Consider the power system described in Definition~\ref{PowerSystem}, where, in addition, we have that $0 < c_{I} < c_{F}$, $c_{F}-c_{F}^{-} = 0$ (free downward regulation), and $F^{-1}(0) = 0$. The stochastic dispatch rule simplifies to:

\vspace{4mm} \begin{center}
{\footnotesize
\begin{tabular}{|l|lll|l|l|}
  \hline
  Rule \# & $p_{W}$ & $p_{I}$ & $p_{F}$ & \multicolumn{2}{c|}{applies if}   \\
  \hline
 7. & $0$ &$0$ & $l$                                                                                             & \multirow{4}{3cm}{$l_3 \leq \mf$} & \multicolumn{1}{l|}{$0 \leq l \leq l_3$} \\
 8. & $0$ &$l - l_3$ & $l_3 $                                                                                   &                                                      & \multicolumn{1}{l|}{$l_3 < l \leq \mi + l_3$} \\
 9. & $0$ & $\mi$ & $l  - \mi$                                                             &                                                      & \multicolumn{1}{l|}{$l_3 + \mi < l \leq \mi + \mf$} \\
10. & $l - \mi -\mf$ & $\mi$ & $\mf$                                        &                                                      & \multicolumn{1}{l|}{$l > \mi + \mf$} \\
\hline
11. & $0$ & $0$ & $l$                                                                                             & \multirow{5}{3cm}{$l_3 > \mf$}      & \multicolumn{1}{l|}{$0 \leq l \leq \mf$} \\
12. & $l - \mf$ & $0$ & $\mf$                                                                     &                                                      &   $\mf \leq l \leq l_4$\\
13. & $l_4 - \mf$ & $l - l_4$ & $\mf$                                                             &                                                      &                        $l_4 \leq l \leq l_4 + \mi$\\
14. & $l - \mf - \mi$ & $\mi$ & $\mf$                                       &                                                      &                       $l > l_4 + \mi$\\
  \hline
\end{tabular}} \end{center}
\vspace{4mm}
\end{corollary}

\proof{Proof.} %
The proof of this corollary directly follows from Theorem~\ref{Prop:ESCF} and the fact that if $c_{F}-c_{F}^{-} = 0$, then $l_1 = F^{-1}(0) = 0 < l_2$.%
\hfill \Halmos \endproof

The case described in this corollary is, perhaps, the most paradigmatic example of a power system in which maximum  efficiency is to be achieved by breaking the merit order in the forward market. Indeed, if the provision of downward regulation does not entail any extra cost to the system (situation that we describe as ``free downward regulation''), the more expensive flexible power capacity should be always dispatched first, even before the cheaper inflexible one, in the hope that the dispatched flexible capacity can be de-committed to accommodate the eventual stochastic power production. One can think of the provision of downward regulation as a sort of arbitrage whereby the flexible generation technology is, in the end, requested not to produce the amount of power that it was scheduled to supply in advance~\citep{Pritchard2010}. Consequently, one can argue about the reasons why the cost of downward regulation should be different from zero, that is, $c_{F}-c_{F}^{-} \neq 0$ (such as the potential existence of nonconvex costs that are not captured in self-commitment based electricity markets).

We now deal with the mirrored case of free upward regulation.

\begin{corollary}[Free upward regulation]\label{Corr:FreeUW}
Consider the power system described in Definition~\ref{PowerSystem}, where, besides, we have that $0 < c_{I} < c_{F}$ and $c_{F}-c_{F}^{+} = 0$ (free upward regulation). The stochastic dispatch rule boils down  to:

\vspace{4mm} \begin{center}
{\footnotesize
\begin{tabular}{|l|lll|l|l|l|}
  \hline
  Rule \# & $p_{W}$ & $p_{I}$ & $p_{F}$ & \multicolumn{3}{c|}{applies if}   \\
  \hline
 1. & $l$ & $0$ & $0$                                                              & \multicolumn{3}{l|}{$0 \leq l \leq l_2$} \\
 2. & $l_2$ &$l - l_2$ & $0$                                                       & \multicolumn{3}{l|}{$l_2 < l \leq \mi + l_2$}\\
 3. & $l-\mi$ & $\mi$ & $0$                                  & \multicolumn{3}{l|}{$\mi + l_2 < l \leq \mi + \mw$}\\
 4. & $\mw$ & $\mi$ &$l - \mw -\mi$                          & \multicolumn{3}{l|}{$\mi + \mw < l \leq \mf + \mi + \mw$}\\
 5. & $l - \mf - \mi$ &$\mi$& $\mf$    & \multicolumn{3}{l|}{$l > \mf + \mi + \mw$}\\
\hline
\end{tabular}} \end{center}
\vspace{4mm}
\end{corollary}

\proof{Proof.} %
Again the proof of this corollary directly follows from Theorem~\ref{Prop:ESCF}, by noting that $c_{F} > c_{I}$ and $c_{F}-c_{F}^{+} = 0$ implies $l_1 = F^{-1}(1) = \mw$. Since $l_2 < \mw$, then $l_2 < l_1$.
\hfill \Halmos \endproof

In the case that the provision of upward regulation does not impose any extra cost on the system (situation that we refer to as ``free upward regulation''), the stochastic market-clearing mechanism~\eqref{Eff} prompts forward dispatch solutions that respect the merit order in the forward market (that is, the more expensive flexible power capacity is only dispatched differently from zero for levels of load at which the inflexible generation technology is requested to produce at maximum capacity). As we will show later, this family of dispatch solutions (namely, those honouring the merit order) can be reproduced, at least partially, by enhanced forms of the conventional two-stage market.

Corollaries~\ref{Corr:FreeDW} and \ref{Corr:FreeUW} define two extreme cases of a broader family of power systems characterized by asymmetric costs for balancing power. In a power system where the cost differential for downward regulation ($c_{F}-c_{F}^{-}$) is sufficiently lower than the cost differential for upward regulation ($c_{F}^{+} - c_{F}$), the conventional two-stage market will prove to be inefficient, because the most economical way of operating the power system is that for which the dispatch of the flexible generating capacity is prioritized over the commitment of the cheaper inflexible one. In the opposite case, where $c_{F}-c_{F}^{-}$ is enough greater than $c_{F}^{+} - c_{F}$, the most cost-efficient forward dispatch complies with the merit order and, therefore, could still be induced by the conventional two-stage market.
Nevertheless, one could argue for the case $c_{F}^{+} - c_{F} > c_{F}-c_{F}^{-}$  to be the general rule, insofar as the provision of upward regulation involves generating more energy than scheduled, while the provision of downward regulation entails not to honour a forward contract.

We now complement Corollaries~\ref{Corr:FreeDW} and \ref{Corr:FreeUW} with the following result, which pertains to the case for which the marginal costs of upward and downward regulation are not marked up with respect to the forward production costs, that is, $c_{F}^{+} - c_{F} = c_{F}-c_{F}^{-} = 0$.

\begin{corollary}[Free upward and downward regulation]\label{Corr:FreeUW-DW}
Consider the power system described in Definition~\ref{PowerSystem}, where, besides, we have that $0 < c_{I} < c_{F}$ and $c_{F}-c_{F}^{+} = c_{F}-c_{F}^{-} = 0$ (free upward and downward regulation). The stochastic dispatch rule reduces  to:

\vspace{4mm}
\begin{center}
 {\footnotesize
\begin{tabular}{|l|ll|l|l|l|}
  \hline
  Rule \# & $p_{W}$ + $p_{F}$ & $p_{I}$ & \multicolumn{3}{c|}{applies if}   \\
  \hline
 1. & $l$ & $0$ &                                                               \multicolumn{3}{l|}{$0 \leq l \leq l_2$} \\
 2. & $l_2$ &$l - l_2$ &                                                         \multicolumn{3}{l|}{$l_2 < l \leq \mi + l_2$}\\
 3. & $l-\mi$ & $\mi$ &                                  \multicolumn{3}{l|}{$\mi + l_2 < l$}\\
\hline
\end{tabular}
}
with $p_{W} \geq 0$ and $p_{F} \geq 0$.
\end{center}
\end{corollary}

\vspace{4mm}
\proof{Proof.} %
This result can be proved by noting that if $c_{F}-c_{F}^{+} = c_{F}-c_{F}^{-} = 0$, then $l_2 = l_3 = l_4$.%
\hfill \Halmos \endproof

Corollary~\eqref{Corr:FreeUW-DW} describes a power system for which dispatch and re-dispatch actions are equally costly. In this case, there is always a forward dispatch solution that satisfies the merit order (by just setting $p_{F} = 0$ in rules 1 and 2). As we will see later, this makes it possible to analyze enhanced variants of the conventional two-stage market that close the cost-efficiency gap with respect to the ideal stochastic market-clearing mechanism~\eqref{Eff}.

In the stochastic dispatch rule given by Theorem~\ref{Prop:ESCF}, the capacities of the flexible and inflexible generation technologies play roles that are as important as their marginal generation costs. In this vein, Corollary~\ref{Corr:CapacityAdequate} below provides two relevant results. The first one pertains to the case of a power system where there is enough flexible power capacity to cover any potential lack of stochastic power production, that is, $\mf \geq \mw$. In an abuse of terminology, we will refer to this instance as ``capacity adequate power system''. The second result constitutes a further simplification of the stochastic dispatch rule that is possible when, in addition, the range of system loads does not exceed the capacity of the inflexible generation technology, that is, $l \leq \mi$.

\begin{corollary}[Capacity adequate power system]\label{Corr:CapacityAdequate}

Consider the power system described in Definition~\ref{PowerSystem}, where, in addition, we have that $c_{I} < c_{F}$ and $\mf \geq \mw$. The stochastic dispatch rule simplifies to:

\vspace{4mm} \begin{center}
{\footnotesize
\begin{tabular}{|l|lll|l|l|l|}
  \hline
  Rule \# & $p_{W}$ & $p_{I}$ & $p_{F}$ & \multicolumn{3}{c|}{applies if}   \\
  \hline
 1. & $l$ & $0$ & $0$                                                             & \multirow{5}{1.2cm}{$l_1 \geq l_2$}  & \multicolumn{2}{l|}{$0 \leq l \leq l_2$} \\
 2. & $l_2$ &$l - l_2$ & $0$                                                       &                                      & \multicolumn{2}{l|}{$l_2 < l \leq \mi + l_2$}\\
 3. & $l-\mi$ & $\mi$ & $0$                                 &                                      & \multicolumn{2}{l|}{$\mi + l_2 < l \leq \mi + l_1$}\\
 4. & $l_1$ & $\mi$ &$l - l_1 -\mi$                        &                                      & \multicolumn{2}{l|}{$\mi + l_1 < l \leq \mf + \mi + l_1$}\\
 5. & $l - \mf - \mi$ &$\mi$& $\mf$ &                                      & \multicolumn{2}{l|}{$l > \mf + \mi + l_1$}\\
\hline
 6. & $l$ & $0$ & $0$                                                             & \multirow{5}{1.2cm}{$l_1 < l_2$}     & \multicolumn{2}{l|}{$0 \leq l \leq l_1$}                                  \\
 7. & $l_1$ &$0$ & $l-l_1$                                                       &                                      & \multicolumn{2}{l|}{$l_1 \leq l \leq l_3$} \\
 8. & $l_1$ &$l - l_3$ & $l_3 - l_1$                                             &                                      & \multicolumn{2}{l|}{$l_3 < l \leq \mi + l_3$} \\
 9. & $l_1$ & $\mi$ & $l - l_1 - \mi$                       &                                     & \multicolumn{2}{l|}{$l_3 + \mi < l \leq \mi + \mf + l_1$} \\
10. & $l - \mi -\mf$ & $\mi$ & $\mf$  &                                     & \multicolumn{2}{l|}{$l > \mi + \mf + l_1$} \\
  \hline
\end{tabular}} \end{center}
\vspace{4mm}
where:
\begin{equation*}
  l_2 = F^{-1}\left(\frac{c_{I}}{c_{F}^{+}}\right) \enskip \text{and} \enskip  l_3 = max\left(l_1, F^{-1}\left(1-\frac{c_{F}-c_{I}}{c_{F}^{-}}\right)\right) \enskip \text{with } l_1 \leq l_3 \leq \mf.
\end{equation*}

If, in addition, we have that $l \leq \mi$, the stochastic dispatch rule further reduces to:

\vspace{4mm} \begin{center}
{\footnotesize
\begin{tabular}{|l|lll|l|l|l|}
  \hline
  Rule \# & $p_{W}$ & $p_{I}$ & $p_{F}$ & \multicolumn{3}{c|}{applies if}   \\
  \hline
 1. & $l$ & $0$ & $0$                                                             & \multirow{2}{1.2cm}{$l_1 \geq l_2$}  & \multicolumn{2}{l|}{$0 \leq l \leq l_2$} \\
 2. & $l_2$ &$l - l_2$ & $0$                                                       &                                      & \multicolumn{2}{l|}{$l_2 < l$}\\
\hline
 6. & $l$ & $0$ & $0$                                                             & \multirow{3}{1.2cm}{$l_1 < l_2$}     & \multicolumn{2}{l|}{$0 \leq l \leq l_1$}                                  \\
 7. & $l_1$ &$0$ & $l-l_1$                                                       &                                      & \multicolumn{2}{l|}{$l_1 \leq l \leq l_3$} \\
 8. & $l_1$ &$l - l_3$ & $l_3 - l_1$                                             &                                      & \multicolumn{2}{l|}{$l_3 < l$} \\
  \hline
\end{tabular}}\end{center}

\end{corollary}

\vspace{4mm} The proof of Corollary \ref{Corr:CapacityAdequate} is provided in Appendix \ref{Proof:CapacityAdequate}.

%Let us assume that $l_3=\mf\gg\mw$. Then, by replacing $l_3$ in \eqref{L3} we have that $c_F + \left(v-c_F^+\right)F\left(\mf-\mf\right) - c_F^- \left(1-F(\mf)\right) \geq c_I$ is satisfied since $\cf>\ci$, $\left(v-c_F^+\right)F\left(0\right) \geq 0$ and  $F(\mf)=1$. This implies that $l_3 \leq \mf$ and thus $l_3 \leq l_1 + \mf$. Therefore, if $l_1 < l_2$ rules 7-10 apply.

For the sake of illustration, later on we will discuss examples for which the conditions of Corollary \ref{Corr:CapacityAdequate} hold. Notice that, despite the fact that these conditions are quite restrictive, they prompt a simplified stochastic dispatch rule that still include cases for which the merit order is broken (see rules 7 and 8).

To conclude this section, we would like to point out the case of symmetric balancing costs, that is, $c_{F}-c_{F}^{+} = c_{F}-c_{F}^{-} = \Delta$. By Equation~\eqref{L1}, we arrive at $l_{1} = F^{-1}(0.5)$, which is the median of the probability distribution that characterizes the stochastic power production. If, besides, this probability distribution is symmetric, $l_1$ is equal to its expected value. Under symmetric balancing costs and $\mf \geq \mw$, the stochastic dispatch rule yields solutions that violate the merit order if $c_{I} > 0.5 c_F^+$ for the case in which $c_{F} > c_{I}$.

\section{A Price-consistent Conventional Two-stage Market (ConvM-VB): The Role of Virtual Bidding} \label{VirtualBidding}

The stochastic two-stage market StoM ensures maximum cost-efficiency \emph{and} price-consistency. In contrast, it may violate the merit order and as a result, dispatch conventional generating units in a loss-making position in the forward market \citep{morales2014electricity}.
In this section, we analyze an enhanced form of the conventional two-stage market~\eqref{Ineff}--\eqref{IneffBM} that, by construction, ensures price consistency and preserves the merit order. Mathematically, we can get price-consistent solutions out of \eqref{Ineff}--\eqref{IneffBM} by freeing variable $p_{W} \geq 0$, which represents the forward dispatch of the stochastic power capacity, while at the same time enforcing the additional constraint $\lambda^f= E_{\Omega}\left[\lambda^{b}(\omega)\right]$. As we show below, to do so, we need to pose the price-consistent conventional two-stage market as a complementarity problem.

In practice, price-consistent solutions can be theoretically obtained from a conventional two-stage market where virtual bidding is allowed and exercised by a risk-neutral arbitrager that has perfect knowledge of the market price distribution (induced by the uncertain power supply). For this purpose, we solve the system of non-linear equations that results from concatenating the KKT conditions of the following optimization problems:

\begin{description}
  \item[Clearing of the forward market]
\begin{subequations}\label{Clearing_DA}
\begin{align}
&\underset{p_{F}, p_{I}, p_{W}}{{\rm  Minimize}} \quad c_F p_{F} + c_I p_{I}\label{OF_Clearing_DA}\\
&{\rm s.t.} \quad  p_{F} + p_{I} + p_W + p_V = l: \lambda^{f}\label{Balance_Clearing_DA}\\
& \phantom{{\rm s.t.}} \quad 0 \leq p_W \leq \uw\\
& \phantom{{\rm s.t.}} \quad 0 \leq p_F \leq \mf\\
& \phantom{{\rm s.t.}} \quad 0 \leq p_I \leq \mi
%
%& \phantom{{\rm s.t.}} \quad p_{F}, p_{W}, p_{I} \geq 0,
%
\end{align}
\end{subequations}
where $p_V$ is the optimal virtual bid coming from the arbitrager's problem below. Notice that we assume that the bid price of the virtual bidder is zero.

  \item[Clearing of the balancing market]
\begin{subequations}\label{Clearing_BM}
\begin{align}
&\underset{p_{F}^+(\omega), p_{F}^-(\omega),  \Delta p_{W}(\omega), \sw}{{\rm  Minimize}} \quad p_{F}^{+}(\omega) c_F^{+} - p_{F}^{-}(\omega) c_F^{-} + v \sw\label{OF_Clearing_BM}\\
& {\rm s.t.} \quad p_{F}^{+}(\omega) - p_{F}^{-}(\omega) + \Delta p_{W}(\omega) + \iv + \sw = 0: \lambda^{b}(\omega)\label{Balance_Clearing_BM}\\
& \phantom{{\rm s.t.}} \quad 0 \leq p_{F}^{+}(\omega) \leq \mf - p_{F}\\
& \phantom{{\rm s.t.}} \quad 0 \leq p_{F}^{-}(\omega) \leq p_{F}\\
& \phantom{{\rm s.t.}} \quad 0 \leq p_{W} + \Delta p_{W}(\omega) \leq W(\omega)\\
& \phantom{{\rm s.t.}} \quad 0 \leq \sw \leq l
%
%& \phantom{{\rm s.t.}} \quad p_{F}^{+}(\omega), p_{F}^{-}(\omega), \sw \geq 0,
\end{align}
\end{subequations}
where $p_{W}$ and $p_{F}$ are given by problem~\eqref{Clearing_DA} and $\Delta p_{V}(\omega)$ is the amount of electricity resold, if positive, or repurchased, if negative, by the virtual bidder in the balancing market.

%\begin{subequations}\label{Clearing_BM2}
%%
%\begin{align}
%%
%&\underset{p^+_{\omega}, p^-_{\omega},  \Delta p^{W}_{\omega}}{{\rm  Minimize}} \quad p^+_{\omega} c^{+} - p^-_{\omega} c^{-}\label{OF_Clearing_BM}\\
%%
%& {\rm s.t.} \quad p^+_{\omega} - p^-_{\omega} + \Delta p^{W}_{\omega} + \Delta p_V(W) = 0: \lambda^{b}_{\omega}\label{Balance_Clearing_BM}\\
%%
%& \phantom{{\rm s.t.}} \quad p^-_{\omega} \leq p_{F}^{\star}\\
%%
%& \phantom{{\rm s.t.}} \quad 0 \leq p_{W}^{\star} + \Delta p^{W}_{\omega} \leq W_{\omega}\\
%%
%& \phantom{{\rm s.t.}} \quad p^+_{\omega}, p^-_{\omega} \geq 0,
%\end{align}
%\end{subequations}

  \item[Arbitrager's problem]
\begin{subequations}\label{Arbitrager}
\begin{align}
&\underset{p_V, \iv}{{\rm  Maximize}} \quad  p_V \lambda^{f} + \int_{\Omega}{\iv \lambda^{b}(\omega)f(\omega) d\omega}\label{OF_Arbitrager}\\
&{\rm s.t.} \quad  p_V + \iv = 0,
\end{align}
\end{subequations}
where $\lambda^{f}$ is given as the Lagrange multiplier associated with constraint~\eqref{Balance_Clearing_DA} and $\lambda^{b}(\omega)$ as the Lagrange multiplier associated with constraint~\eqref{Balance_Clearing_BM}.
\end{description}

We now characterize the solution to the system of KKT conditions associated with the convex optimization problems~\eqref{Clearing_DA}, \eqref{Clearing_BM} and \eqref{Arbitrager} (the so-called \emph{short-run equilibrium solution}). Notice that the forward dispatch associated with this solution exhibits a fundamental property, namely, it satisfies the merit order, while leading to price-consistency.

\begin{theorem}[A price-consistent conventional two-stage market]\label{Prop:VB}
Consider the stylized power system described in Definition~\ref{PowerSystem}. Consider also the equilibrium problem that results from simultaneously enforcing the optimality conditions of problems~\eqref{Clearing_DA}, \eqref{Clearing_BM} and \eqref{Arbitrager}, where $c_{I} < c_{F}$. Define the constants
\begin{align}
l_2 &:= \underset{l \geq 0}{\rm min} \enskip l: \left(v-c_F^+\right)F\left(l-\mf\right) + c_F^+F(l) \geq c_I\\
l_5 &:= \underset{l \geq 0}{\rm min} \enskip l: \left(v-c_F^+\right)F\left(l-\mf- \mi\right) + c_F^+F(l- \mi) \geq c_F \label{L5}\\
l_6 &:= \underset{l \geq 0}{\rm min} \enskip l: \left(v-c_F^-\right)F\left(l-\mf- \mi\right) + c_F^-F(l- \mi) \geq c_F \label{L6}
\end{align}
and the function $l_7: \mathbb{R}_{\geq 0} \rightarrow \mathbb{R}_{\geq 0}$,
%
%\begin{equation}\label{VBsol}
%  \varphi(l) := \left\{y: \left(v-c_{F}^{+}\right) F\left(l-\mi-\mf\right) + \left(c_{F}^{+}-c_{F}^{-}\right)F\left(y\right) + c_{F}^{-} F\left(l-\mi\right) = c_F\right\}.
%\end{equation}

\begin{equation}\label{VBsol}
l_7(l) = \Fi{\dfrac{\cf - \pp{v-\cfu}\F{l-\mi-\mf} - \cfd\F{l-\mi}}{\cfu - \cfd}}.
\end{equation}
%

%\textcolor{red}{You have to delete one of the two previous equations. I may like better $l_7$, but I guess it's not too important.}

Note that, by construction, $l_2 \leq l_5 \leq l_6$. The equilibrium solution is given by:

\vspace{4mm} \begin{center}
{\footnotesize
\begin{tabular}{|l|lll|l|}
\hline
Rule \# & $\pw + \pv$ & $\pii$ & $\pf$ & applies if \\
\hline
1. & $l$ & $0$ & $0$ &  $0 \leq l \leq l_2$ \\
2. & $l_2$ & $l-l_2$ & $0$ &  $l_2 < l \leq l_2+\mi$ \\
3. & $l-\mi$ & $\mi$ & $0$ & $l_2 + \mi < l \leq l_5$ \\
4. & $l_7(l)$ & $\mi$ & $l - l_7(l) - \mi$ & $l_5 < l \leq l_6$ \\
5. & $l-\mf-\mi$ & $\mi$ & $\mf$  & $l_6 < l$ \\
\hline
\end{tabular}
}
in all cases with $p_{W} \leq \uw$. \end{center}
\end{theorem}

\vspace{4mm}
A proof for this theorem can be found in Appendix~\ref{Proof:VB}.

We can now use this theorem to provide conditions under which the conventional two-stage market with virtual bidding results in maximum cost-efficiency.

\begin{corollary}[Merit-order dispatch solution with virtual bidding]\label{Corr:MO}
Consider the power system described in Definition~\ref{PowerSystem}, where, in addition, we have that $\ci < \cf$. The expected system operating cost associated with the \emph{price-consistent merit-order} dispatch solution is equal to that of the stochastic dispatch solution in any of the following cases:

\begin{enumerate}
  \item $l_1 \geq l_2$ and $0 \leq l \leq min\left(\mi + l_1, l_5\right)$;
  \item $\nvs \geq \frac{c_I}{c_{F}^{+}}$ and $0 \leq l \leq \mi + l_1$ provided that $\mf > \mw$ (capacity adequate power system).
  \item $l \notin \left(l_5, \mw + \mi + \mf\right)$ if $\cfu = \cf$ (free upward regulation).
  \item $0 \leq l \leq \mi+\mw$ or $l \geq \mw + \mi + \mf$, if $\cfu = \cf$ and $\mf > \mw$  (i.e., free upward regulation in a capacity adequate power system).
  \item $\cf=\cfu=\cfd$ (free upward and downward regulation).
\end{enumerate}

Furthermore, if $\cf=\cfd$ (free downward regulation), $\cfu - \cf > 0$, and $\Fi{0} = 0$, then the price-consistent conventional two-stage market does not deliver maximum cost-efficiency except for $l \geq max(l_{6}, l_4 + \mi)$ when $l_3 > \mf$, and $l \geq max(l_{6}, \mf + \mi)$ otherwise.
% If this system is such that $\frac{c_I}{c_{F}^{+}} \leq \frac{c_F-c_{F}^{-}}{c_F^{+}-c_{F}^{-}}$, then the dispatch solution induced by the merit order with virtual bidding for the inflexible and flexible power capacity and, hence, the resulting expected system operating cost are the same as those given by the efficient market~\eqref{Eff}.
\end{corollary}

The proof of Corollary \ref{Corr:MO} is provided in Appendix \ref{Proof:MO}.

It is apparent that enforcing price consistency, through the introduction of virtual bidding, makes the dispatch solution of the conventional two-stage market substantially more intricate. Most importantly, virtual bidding enhances the efficiency of the conventional two-stage market to such an extent that the cost-efficiency gap with respect to the stochastic dispatch solution is closed under certain conditions. In essence, the ability of virtual bidding to close this gap is mostly determined by the relation between the characteristic system constants $l_1$ and $l_2$. On the one hand, the value of $l_1$ (which is the solution to a news-vendor-type of problem as remarked in Appendix~\ref{Proof:ESCF}) is driven by the relation between the cost of providing downward regulation versus the cost of providing upward regulation. On the other, the value of $l_2$ weighs the cost of dispatching the inflexible power capacity against the cost of dispatching the stochastic power capacity. The former action incurs a marginal cost of $\ci$, while the latter entails a probable marginal cost of $\cfu$, or $v$ in the worst-case scenario, given that any potential shortage of stochastic power production will have to be covered with upward regulation from the flexible generation technology and/or load curtailment.

Thus, the effectiveness of virtual bidding to increase the cost-efficiency of the conventional two-stage market is greater for those power systems with a fleet of cheap inflexible power plants and comparatively costly means for downward regulation ($l_1 \geq l_2$). In such a case, virtual bidding can even drive the conventional two-stage market to maximum cost-efficiency for a range of system loads.

As previously mentioned, the conventional two-stage market with virtual bidding yields dispatch solutions that respect the merit order, while resulting in price-consistency. This leads to an interesting conclusion. Even under those conditions for which the merit order is preserved in the stochastic dispatch rule ($l_1 \geq l_2$), a price-consistent conventional two-stage market does not necessarily deliver maximum cost-efficiency. This happens, for example, in the load range $l_5 < l < l_6$. This implies that a price-consistent merit-order market solution is \emph{not}, in general, the merit-order market solution that minimizes the expected system operating cost. This will become more clear in the illustrative example below and when we introduce, later on, the conventional two-stage market with centralized dispatch of stochastic power production.

\subsection{Example 1}\label{Ex1}

In this example we consider the stylized power system described in
Definition~\ref{PowerSystem} with $v = \$1000$/MWh.
%
%Consider a power system with infinite transmission capacity. The generation portfolio of the system includes wind power production with a limited capacity $\mw$ of 100 MW, and flexible and inflexible power generation technologies, both with infinite capacity.
%
%The marginal cost of the inflexible and flexible power production is $c_{I} = \$30$/MWh and  $c_{F} = \$35$/MWh, respectively. Furthermore, the flexible power capacity is able to provide upward and downward regulation at a cost of $c_{F}^{+} = \$40$/MWh and  $c_{F}^{-} = \$30$/MWh, in that order.
%
The uncertain power supply comes from a wind power farm whose
capacity factor is assumed to follow a Beta distribution. The mean
($\kappa$) and standard deviation ($\sigma$) of this distribution
are linked together through the empirical relationship
\eqref{StandDeviation} provided in~\citet{fabbri2005assessment}.
The shape parameters $\alpha$ and $\beta$ of the Beta distribution
modeling the wind power capacity factor are, consequently,
computed according to \eqref{Beta}.
\begin{equation}
\sigma = 0.01837 + 0.20355 \cdot \kappa\;.
\label{StandDeviation}
\end{equation}
\begin{equation}\label{Beta}
\alpha = \frac{\left(1-\kappa\right)\cdot \kappa\cdot \kappa}{\sigma^2}-\kappa, \quad \quad \quad \quad \beta = \alpha \left( \frac{1-\kappa}{\kappa}   \right).
\end{equation}

The expected wind power production is hence given by the product of the predicted wind power capacity factor $\kappa$ and the wind power capacity $\mw$.

We now examine and compare the stochastic two-stage market and the
conventional two-stage market with and without virtual bidding. For
this purpose, we consider the three illustrative cases collated in Table \ref{Tab:data1} and discussed below.
Table~\ref{Tab:results1} provides the dispatch solution
$\left(\pii, \pf, \pw\right)$, or $\left(\pii, \pf, \pw
+\pv\right)$, and the expected total operating cost prompted by
the efficient two-stage market and the conventional two-stage
market with and without virtual bidding. This table also includes
the forward price $\lf$ and the expected real-time price
$\lbh$ under each of the markets, and the incremental
running cost in percentage with respect to the most cost-efficient
dispatch solution, which is the one provided by the stochastic
two-stage market, logically. In the sequel, we will refer to this
incremental cost simply as \emph{efficiency gap}.

\begin{table} \setlength{\tabcolsep}{8pt} \centering
\begin{tabular}{cccccccccc}
\hline
& $\mi$ & $\mf$ & $\mw$ & $\ci$ & $\cf$ & $\cfu$ & $\cfd$ & $l$ & $\kappa$ \\
\hline
Case a) & 500 & 500 & 100 & 30 & 35 & 35 & 30 & 250 & 0.5 \\
Case b) & 500 & 500 & 100 & 30 & 35 & 40 & 30 & 250 & 0.5 \\
Case c) & 100 & 50  & 100 & 30 & 35 & 35 & 30 & 170 & 0.5 \\
\hline
\end{tabular} \caption{Example 1: Data. System load and capacities are given in MW and marginal costs in \$/MWh.} \label{Tab:data1}
\end{table}

\begin{itemize}
\item Case a) considers a power system with asymmetric balancing costs,
in which upward regulation is less costly than downward
regulation. More specifically, this instance corresponds to a
capacity adequate power system with free upward regulation. Hence,
the results in Corollary~\ref{Corr:FreeUW} apply and, accordingly,
the stochastic two-stage market provides a solution that respects
the merit order, because the capacity of the more expensive
flexible power generation is reserved to deploy upward regulation
in real time, when needed. Naturally, the conventional two-stage
market always yields, by construction, a dispatch solution that
satisfies the merit order. However, in the case that virtual
bidding is not allowed, the conventional two-stage market clears
the expected wind power production, which is, in general, a
suboptimal dispatch decision in terms of market efficiency and as
such, results in an efficiency gap greater than zero. On the
contrary, the introduction of virtual bidding not only induces a
price-consistent market solution, but also fully closes the
efficiency gap. This result is consistent with Claim 4 in
Corollary~\ref{Corr:MO}, since $0 \leq l \leq \mi + \mw$, with $l
= 250$ MW, $\mi = 500$ MW, and $\mw = 100$ MW.

\item Case b) represents a capacity adequate power system identical to that of Case a),
except for the fact that the provision of upward regulation is now
costly, that is, $\cfu > \cf$. In particular, we have that
$\cfu-\cf = \cf - \cfd =$ \$5/MWh (symmetric balancing costs).
Besides, since $l < \mi$, the second result in Corollary~\ref{Corr:CapacityAdequate} applies.
Lastly, because $\ci > 0.5\cfu$, the stochastic two-stage market
provides a dispatch solution that breaks the merit order.
Evidently, the conventional two-stage market cannot replicate such
a dispatch, irrespective of whether virtual bidding is permitted
or not. Nevertheless, virtual bidding reduces the efficiency gap,
while ensuring price consistency.
\item Case c) is also similar to Case a), but with reduced generating capacities and system
load. In fact, this instance corresponds to a capacity
\emph{inadequate} power system with free upward regulation.
Consequently, results from Corollary~\ref{Corr:FreeUW} hold and
the stochastic two-stage market prompts a dispatch solution that
does comply with the merit order. However, as opposed to what happens
in Case a), not even virtual bidding is able to close the
efficiency gap under the conventional two-stage market in this
instance. This shows that a price-consistent merit-order-based
market does not necessarily deliver maximum cost-efficiency.
\end{itemize}

\begin{table} \setlength{\tabcolsep}{8pt} \centering
\begin{tabular}{ccccccccc}
\hline
& & $\pw (+ \pv)$ & $\pii$ & $\pf$ & $\lf$ & $\lbh$ & cost & gap \\
\hline
 & {StoM} & 63.5 & 186.5 &  0 & 30 & 30 & 6095 & - \\
Case a) & {ConvM} & 50 & 200 & 0 & 30 & 17.5  & 6170 & 1.2\% \\
 & {ConvM-VB} & 63.5 & 186.5 & 0 & 30 & 30 & 6095 & 0\% \\
\hline
 & {StoM} & 50.5 & 188.5 & 11 & 30 & 30 &  6139 & - \\
Case b) & {ConvM} & 50 & 200 & 0 & 30 & 20 & 6195 & 0.9\% \\
 & {ConvM-VB} & 58.5 & 191.5 & 0 & 30 & 30 & 6154 & 0.2\% \\
\hline
 & {StoM} & 70 & 100 &  0 & 37.09 & 37.09 & 3717 & -\\
Case c) & {ConvM} & 50 & 100 & 20 & 35 & 34.83 & 3740 & 0.6\% \\
 & {ConvM-VB} & 51.5 & 100 & 18.5 & 35 & 35 & 3737 & 0.5\% \\
\hline
\end{tabular} \caption{Example 1: Results. Power dispatch values are given in MW, prices in \$/MWH, cost in \$/h, and incremental cost (efficiency gap) in percentage.} \label{Tab:results1}
\end{table}

To conclude this example,  we investigate how the mean capacity
factor $\kappa$ and the forecast horizon impacts the efficiency
gap. To this aim, we consider Case d), for which data is provided
in Table \ref{Tab:data2}. Note that this case corresponds to a capacity adequate power system with asymmetric balancing costs. To be more precise, the provision of downward regulation in this system does not entail any extra operating cost since $\cfd = \cf =$ \$35/MWh. Therefore, results from Corollary~\ref{Corr:FreeDW}, which in general render a dispatch solution that violates the merit order,
apply. Accordingly, we should expect that the dispatch solutions
induced by the conventional two-stage market, both with and
without virtual bidding, feature a nonzero efficiency gap in this
case.

\begin{table} \setlength{\tabcolsep}{8pt} \centering
\begin{tabular}{ccccccccc}
\hline
& $\mi$ & $\mf$ & $\mw$ & $\ci$ & $\cf$ & $\cfu$ & $\cfd$ & $l$ \\
\hline
Case d) & 500 & 500 & 100 & 30 & 35 & 40 & 35 & 250 \\
\hline
\end{tabular} \caption{Example 1: Data for Case d).} \label{Tab:data2}
\end{table}

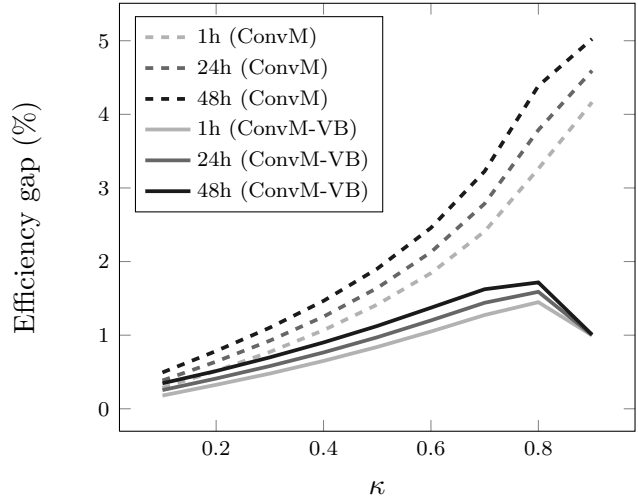
\begin{figure}
\centering
\begin{tikzpicture}[scale=1]
    \begin{axis}[
    ticklabel style = {font=\scriptsize},
    xlabel = $\kappa$,
    ylabel = Efficiency gap (\%),
    legend cell align=left,
    legend style={legend pos=north west,font=\scriptsize},
    legend entries={1h (ConvM),24h (ConvM),48h (ConvM),1h (ConvM-VB),24h (ConvM-VB),48h (ConvM-VB)}]
    \addplot[color=black!30,line width=0.5mm, dashed] table [x="fore", y="inef"] {table4_1.tex};
    \addplot[color=black!60,line width=0.5mm, dashed] table [x="fore", y="inef"] {table4_24.tex};
    \addplot[color=black!90,line width=0.5mm, dashed] table [x="fore", y="inef"] {table4_48.tex};
    \addplot[color=black!30,line width=0.5mm] table [x="fore", y="vb"] {table4_1.tex};
    \addplot[color=black!60,line width=0.5mm] table [x="fore", y="vb"] {table4_24.tex};
    \addplot[color=black!90,line width=0.5mm] table [x="fore", y="vb"] {table4_48.tex};
    \end{axis}
\end{tikzpicture} \caption{Efficiency gap for different expected wind power capacity factors and forecast horizons} \label{Fig:EfficGap}
\end{figure}

This is confirmed by the plots in Figure \ref{Fig:EfficGap}, which
illustrate such an efficiency gap in percentage (on the y-axis)
for different values of $\kappa$ (on the x-axis) and forecast
horizons (i.e., time spans in between the clearings of the forward
and the real-time markets), namely, 1, 24 and 48 hours. Each
forecast horizon is associated with a different empirical
relationship between $\sigma$ and $\kappa$ in the form
of~\eqref{Beta}, all of which have been taken
from~\citet{fabbri2005assessment}. The continuous and dashed lines
correspond to the cases with and without virtual bidding, respectively.

Two observations are in order: First, the efficiency gap increases
as the time distance in between the closures of the forward and
the real-time market augments. This is  hardly a surprising result
related to the fact that uncertainty in wind power production
grows as does the forecast lead time. Second, virtual bidding, or
more generally, a price-consistent conventional two-stage market,
can substantially reduce the efficiency gap as compared to the
case of a conventional market that does not ensure price
consistency, particularly in those situations where the
contribution of the stochastic power production is expected to be
important.

Finally, it is worth noting that the efficiency gap decreases for
high values of $\kappa$ in the case where virtual bidding is
allowed. This is a direct consequence of the fact that, as the
value of $\kappa$ grows, the resulting wind power distribution
becomes more and more skewed towards high wind power production
values. Intuitively speaking, the probability mass concentrates
more and more towards values of wind power production closer to the wind farm capacity
$\mw$ as $\kappa$ increases, a phenomenon that virtual bidding
captures and takes advantage of.

\section{A Conventional Two-stage Market with Centralized Dispatch of the Stochastic Power Production (ConvM-CD)} \label{Centralized}
In the sequel, we analyze a variant of the conventional two-stage market~\eqref{Ineff}--\eqref{IneffBM}, whereby the forward dispatch $\pw$ of stochastic power capacity is centralized and determined by a non-profit and all-knowing entity that seeks to minimize the expected total system operating cost. A transmission system operator, for example, could take on this role. The mathematical model that we present next to simulate this market organization is inspired from the one introduced in~\cite{morales2014electricity} and takes the form of the following bilevel linear programming problem.
\begin{subequations}\label{EJOR}
\begin{align}
& \hspace{0mm} \underset{\pf,\pii,\pw; \pfu, \pfd, \iw, \sw}{{\rm Minimize}}\quad \ci \pii + \cf \pf + \int_{\Omega}{\left(\cfu \pfu - \cfd \pfd + v \sw\right) f(\omega) d\omega} \label{EJOR_OF}\\
& \hspace{0mm}{\rm s.t.}\;\; p_{F}^{+}(\omega) - p_{F}^{-}(\omega) + \Delta p_{W}(\omega) + \sw = 0: \lambda^{b}(\omega), \enskip \forall \omega \in \Omega \label{EJOR_First}\\
& \hspace{0mm}\phantom{s.t.}\;\; 0 \leq p_{F}^{+}(\omega) \leq \mf - p_{F}, \enskip \forall \omega \in \Omega\\
& \hspace{0mm}\phantom{s.t.}\;\; 0 \leq p_{F}^{-}(\omega) \leq p_{F}, \enskip \forall \omega \in \Omega\\
& \hspace{0mm}\phantom{s.t.}\;\; 0 \leq p_{W}  \leq \mw, \enskip \forall \omega \in \Omega\\
& \hspace{0mm}\phantom{s.t.}\;\; 0 \leq p_{W} + \Delta p_{W}(\omega) \leq W(\omega), \enskip \forall \omega \in \Omega\\
& \hspace{0mm}\phantom{s.t.}\;\; 0 \leq \sw \leq l, \enskip \forall \omega \in \Omega\label{EJOR_Final}
\end{align} \vspace{-12mm}
\begin{empheq}[left=\hspace{-55mm}{\left(\pf, \pii\right) \in \text{arg}}\empheqlbrace \enskip,right=\enskip\empheqrbrace .]{align}
& \underset{x_{F} , x_{I}}{{\rm Minimize}} \quad c_F x_{F} + c_I x_{I} \label{EJOR_LL_OF}\\
& {\rm s.t.} \;\;  x_{F} + x_{I} + p_W = l: \lambda^{f}\\
& \phantom{{\rm s.t.} \;\;} 0 \leq x_F \leq \mf\\
& \phantom{{\rm s.t.} \;\;}0 \leq x_I \leq \mi \label{EJOR_LL_Final}
%
%& \phantom{{\rm s.t.} \;\;}x_{F}, x_{I} \geq 0\label{EJOR_LL_Final}\
%
\end{empheq}
\end{subequations}

Essentially, the lower-level
problem~\eqref{EJOR_LL_OF}--\eqref{EJOR_LL_Final} models the
clearing of the forward market as a function of the amount $p_W$
of stochastic power production that is dispatched. This amount is
determined by the upper-level problem
\eqref{EJOR_OF}--\eqref{EJOR_Final} with a view to minimizing the
expected total system operating cost~\eqref{EJOR_OF}. To this aim,
the upper-level problem explicitly anticipates, through
\eqref{EJOR_First}--\eqref{EJOR_Final}, the projected re-dispatch
actions that will need to be undertaken for any possible
realization $\omega$ of the stochastic power production. What is most important, though, is that the bilevel linear model~\eqref{EJOR} yields, by construction, a market solution that, while preserving the merit order (due to the enforcement of~\eqref{EJOR_LL_OF}--\eqref{EJOR_LL_Final}), results in maximum cost-efficiency.

The following theorem characterizes the solution to the bilevel linear programming problem~\eqref{EJOR}.
\begin{theorem}[A conventional two-stage market with centralized dispatch of $\pw$]\label{Prop:EJOR}
Consider the power system described in Definition~\ref{PowerSystem} with $c_{I} < c_{F}$. The solution to the bilevel linear programming problem~\eqref{EJOR} is given by:

\vspace{4mm}
\begin{center}
{\footnotesize
\begin{tabular}{|c|lll|l|l|}
\cline{1-6}
Rule \# & $\pw$ & $\pii$ & $\pf$ & \multicolumn{2}{c|}{applies if}\\
\cline{1-6}
1. & $l$ & $0$ & $0$ && $0 \leq l \leq l_2$\\
2. & $l_2$ & $l-l_2$ & $0$ && $l_2 \leq l \leq \mi + l_2$\\
3. & $l-\mi$ & $\mi$ & $0$ & $l_1 \geq l_2$ & $\mi + l_2 \leq l \leq \mi + l_1$\\
4. & $l_1$ & $\mi$ & $l-\mi-l_1$ && $\mi + l_1 \leq l \leq \mi+\mf + l_1$\\
5. & $l-\mi-\mf$ & $\mi$ & $\mf$ && $\mi+\mf + l_1 \leq l$\\
\cline{1-6}
6. &$l$ & $0$ & $0$ &\multirow{4}{2cm}{$l_1 < l_2$}& $0 \leq l \leq \min(l_2,l_8)$\\
7. &$l_2$ & $l-l_2$ & $0$ && $l_2 \leq l \leq l_8$\\
8. &$l_1$ & $\mi$ & $l-\mi-l_1$ && $l_8 \leq l \leq \mi+\mf + l_1$\\
9. &$l-\mi-\mf$ & $\mi$ & $\mf$ & & $\max(l_8,\mi+\mf+l_1) \leq l$\\
\cline{1-6}
\end{tabular}}
\end{center}
\vspace{4mm}

where constant $l_8$ is defined when $l_1 < l_2$ as:
\begin{align}
& l_8 := x: l_1 + \mi \leq x \leq l_2 + \mi \ \text{and} \ \ci\pp{\mi-\max(x-l_2,0)} - \int_{\min(x-\mi,l_1)}^{\min(x-\mi,l_1+\mf)}(A(s)-B(s))ds - \nonumber \\
& \qquad - \int_{\min(x-\mi,l_1+\mf)}^{x-\mi}(A(s)-C(s))ds - \int_{x-\mi}^{\min(x,l_2)}A(s)ds = 0, \\
& \text{with}\\
& \qquad A(s) := (v-\cfu)\F{s-\mf}+\cfu\F{s} \\
& \qquad B(s) := \cf + \pp{v-\cfu}\F{s-\mf}-\cfd\FF{s} \\
& \qquad C(s) := \pp{v-\cfd}\F{s-\mf}+\cfd\F{s}
\end{align}
\end{theorem}

A proof of this theorem is provided in Appendix~\ref{Proof:EJOR}.

As we did for the case of the conventional two-stage market with virtual bidding, we can now use Theorem~\ref{Prop:EJOR} to identify conditions under which a conventional two-stage market with centralized dispatch of the stochastic power production delivers maximum cost-efficiency.

\begin{corollary}[Cost-efficient merit-order dispatch solution]\label{Corr:MO_EJOR}
Consider the power system described in Definition~\ref{PowerSystem}, where, in addition, we have that $\ci < \cf$. The expected system operating cost associated with the \emph{merit-order} dispatch solution provided by the bilevel linear programming problem~\eqref{EJOR} is equal to that of the stochastic dispatch solution in any of the following cases:

\begin{enumerate}
  \item $l_1 \geq l_2$;
    \item Whenever the stochastic dispatch solution satisfies the merit order;
  \item Whenever the stochastic dispatch solution is such that $\pii = \mi$;
  \item $\cf=\cfu=\cfd$ (free upward and downward regulation).
\end{enumerate}
\end{corollary}

%\tred{If $\ci < \cf$, any dispatch with $\pii = \mi$ follows the merit order according to our definition. Therefore, I think 2 and 3 should be merged.}

The proof of Corollary \ref{Corr:MO_EJOR} is provided in Appendix \ref{Proof:MO_EJOR}.

Interestingly, the previous results show that all the stochastic
dispatch solutions that do not violate the merit order could be
recovered, in principle, from a conventional two-stage market, if
we let an all-knowing and social-welfare-maximizer entity (for
example, a TSO) control the dispatch of the stochastic power
production. This implies that a TSO could close the efficiency gap
between the stochastic and the conventional two-stage markets in
many cases where virtual bidding fails to do so. On the other
hand, neither a TSO nor virtual bidding are able to close such
a gap in all those cases where maximum cost-efficiency requires
breaking the merit order. In those cases, though, a TSO would be
able to reduce the efficiency loss further than virtual bidding.
In this vein, when compared with the results provided in Theorems
\ref{Prop:ESCF} and \ref{Prop:EJOR}, the dispatch solution induced
by a conventional two-stage market with virtual bidding reveals
that price-consistency does \emph{not} necessarily implies maximum
cost-efficiency. But it is more interesting yet to notice that
this conclusion also means that the conventional two-stage market with
centralized dispatch of the stochastic power production can
deliver dispatch solutions that are price inconsistent. Indeed,
Theorem \ref{Prop:VB} tells us that this may happen in the range
of loads in between $l_5$ and $l_6$. We illustrate one of these
cases in the example of Section~\ref{Ex2}.

We conclude with the corollary below, which states that the price-consistent and the cost-efficient merit-order dispatch solutions are equivalent when the inflexible generating capacity is sufficiently large.

\begin{corollary}[Merit-order dispatch solution with large inflexible power capacity]\label{Corr:InfCap}
Consider the stylized power system described in Definition~\ref{PowerSystem} with $\mi \geq l$. Then, the solutions of the bilevel model \eqref{EJOR} and the complementarity model consisting of the KKT conditions of problems \eqref{Clearing_DA}, \eqref{Clearing_BM} and \eqref{Arbitrager} are the same. The dispatch rule in that case is given by:

\vspace{4mm}
\begin{center}
{\footnotesize
\begin{tabular}{|l|lll|l|}
\hline
Rule \# & $\pw$ & $\pii$ & $\pf$ & applies if \\
\hline
1. & $l$ & $0$ & $0$ &  $0 \leq l \leq l_2$ \\
2. & $l_2$ & $l-l_2$ & $0$ &  $l_2 < l$ \\
\hline
\end{tabular}
}
\end{center}
\end{corollary}

\proof{Proof.}
From Theorem \ref{Prop:VB}, only rules 1 and 2 apply for $\mi \geq l$. From Theorem \ref{Prop:EJOR} with $l_1 \geq l_2$ and $\mi \geq l$, only rules 1 and 2 apply. Furthermore, if $l_1 < l_2$ and $\mi \geq l$, then $\l_8 \geq l$ and consequently, only rules 6 and 7 in Theorem \ref{Prop:EJOR} apply.
\hfill \Halmos \endproof

\subsection{Example 2}\label{Ex2}

The main purpose of this example is to illustrate the key
differences between the merit-order dispatch solutions provided by
the two enhanced variants of the conventional two-stage market
that we have examined in this paper, namely, that in which virtual
bidding is allowed (ConvM-VB) and that in which the stochastic
power production is centrally dispatched (ConvM-CD). To this aim,
we consider Case e), for which data is provided in Table
\ref{Tab:data3}. Note that this new case is identical to Case d),
but for a different system load of 155 MW. Results for the several
markets analyzed are collated in Table \ref{Tab:results2}. These
results consist of dispatched quantities $(\pw, \pii,\pf)$ or
$(\pw+\pv, \pii,\pf)$, as appropriate; forward price ($\lf$);
expected real-time price ($\lbh$); expected total
system operating cost (cost), and efficiency gap.

Case e) corresponds to a capacity inadequate power system with
asymmetric balancing costs. In particular, the provision of
downward regulation in this system is free in the sense that $\cfd
= \cf$. Therefore, the stochastic dispatch rule is governed by the
results of Corollary~\ref{Corr:FreeDW}, which generally prompt
dispatch solutions that break the merit order, as is the case here
(see row labeled ``StoM'' in  Table~\ref{Tab:results2}).
Hence, none of the enhanced variants of the conventional
two-stage market we consider, i.e., ConvM-VB and ConvM-CD, are able
to close the efficiency gap (see column tagged as ``gap'' in Table~\ref{Tab:results2}).

The stochastic two-stage market delivers minimum expected system
operating cost at the same time that it guarantees price
consistency. However, in order to produce a price-consistent
\emph{and} maximum cost-efficient dispatch solution, it has to
violate the merit order. This has important negative implications
towards the actual implementation of this market-clearing
procedure \citep{morales2014electricity}. As an example, it is
easy to infer from the results in Table~\ref{Tab:results2} that
the flexible power capacity is dispatched in a loss-making
position under the stochastic two-stage market, because the
forward price equals \$30/MWh, whereas the marginal production
cost of the flexible power generation technology is \$35/MWh.

On the other hand, the conventional two-stage market with virtual
bidding (ConvM-VB) produces a price-consistent dispatch solution
that respects the merit order, but that is not the best in terms
of cost efficiency. Indeed, the market solution provided by
ConvM-VB is even less cost-efficient than the one delivered by the
plain conventional two-stage market (ConvM). Although this should
be regarded as a rare case, it shows that ensuring price
consistency does not necessarily lead to higher market efficiency
(understood as the minimization of system operating costs).

Lastly, the conventional market with centralized dispatch of the
stochastic power production (ConvM-CD) yields the most
cost-efficient dispatch solution among those that comply with the
merit order. Note that ConvM-CD gets to reduce the efficiency gap
by half. To do so, however, it must relinquish price-consistency.

In short, preserving the merit order in forward electricity
markets with uncertain supply implies giving up on cost-efficiency or price-consistency (or both).

\begin{table} \setlength{\tabcolsep}{8pt} \centering
\begin{tabular}{cccccccccc}
\hline
& $\mi$ & $\mf$ & $\mw$ & $\ci$ & $\cf$ & $\cfu$ & $\cfd$ & $l$ & $\kappa$ \\
\hline
Case e) & 100 & 50 & 100 & 30 & 35 & 40 & 35 & 155 & 0.5 \\
\hline
\end{tabular} \caption{Example 2: Data. System load and capacities are given in MW and marginal costs in \$/MWh.} \label{Tab:data3}
\end{table}

\begin{table} \setlength{\tabcolsep}{8pt} \centering
\begin{tabular}{ccccccccc}
\hline
& & $\pw$ ($+\pv$) & $\pii$ & $\pf$ & $\lf$ & $\lbh$ & cost & gap \\
\hline
\multirow{4}{1.5cm}{Case e)} & {StoM} & 18.5 & 92.5 & 44 & 30 & 30 & 3245 & - \\
& {ConvM} & 50 & 100 & 5 &  35 & 25.42 & 3296 & 1.6\% \\
 & {ConvM-VB} & 58.5 &  96.5 &  0 &  30 & 30 & 3304 & 1.8\%  \\
 & {ConvM-CD} & 18.5 & 100 & 36.5 & 35 & 22.93 & 3272 & 0.8\% \\

\hline
\end{tabular} \caption{Example 2: Results. Power dispatch values are given in MW, prices in \$/MWh, cost in \$/h, and incremental cost (efficiency gap) in percentage.} \label{Tab:results2}
\end{table}

\section{Conclusions} \label{Conclusion}

The overall message that emerges from our analysis is that the
concepts of \emph{merit order}, \emph{cost-efficiency}, and
\emph{price-consistency} are conflicting requirements that cannot
be met together in an electricity market. Indeed, preserving the
merit order generally comes at the expense of market efficiency,
either in the form of price-inconsistent or cost-inefficient
market solutions.

We have reached this conclusion after examining four types of
two-stage markets with uncertain supply, in which only one or two
of the above requirements can be guaranteed \emph{in general}. We
have provided mathematical formulations for these four types of
markets and explain how these mathematical abstractions translate
into practice.

Our study has also revealed general conditions under which
preserving the merit order is most likely to jeopardize
cost-efficiency and/or price-consistency. Interestingly, this
seems to be the case of power systems that are capacity inadequate
(broadly understood as a power system where the complete lost of
the uncertain supply cannot be covered with economical flexible
power generation) and/or where the provision of upward regulation
is costly, while the provision of downward regulation
comparatively imposes little or no \emph{extra} cost at all to the
system. Likewise, our study has identified general conditions
under which the merit order, price-consistency and cost-efficiency
can indeed be met together. This is, for example, the case of
power systems where the implementation of real-time adjustments do
not entail opportunity costs.

In an attempt to keep our investigation essentially analytical, we
have built our market models on a stylized power system with
infinite transmission capacity. Consequently, a natural avenue for
future research is to elucidate
whether more realistic assumptions on the underlying power system
can diminish the loss of efficiency caused by the merit order or
even limit the cases of inefficiency to a few ``degenerate'' ones.
To this end, we are most likely to abandon our analytical approach
and make use of computational simulation instead.

On a different front, our analysis has considered energy-only
electricity markets. Therefore, another logical direction for
future work is to investigate whether the consideration of other
types of market mechanisms could reduce or even nullify the loss
of efficiency associated with the preservation of the merit order.
In this regard, we conjecture, based on our results, that
appropriate markets for downward operating reserve might do the
trick.

% Appendix here
% Options are (1) APPENDIX (with or without general title) or
%             (2) APPENDICES (if it has more than one unrelated sections)
% Outcomment the appropriate case if necessary
%
% \begin{APPENDIX}{<Title of the Appendix>}
% \end{APPENDIX}
%
%   or
%
\begin{APPENDICES}
\section{Proof of Proposition \ref{Prop:Cost2Stage}} \label{Proof:Cost2Stage}

%\begin{proposition} \label{pro_bal_cost}
%Consider a power system as described in Theorem XX with $0 < \cfd < \cfu < v$. The expected balancing cost for given values of $\pf$ and $\pw$ is equal to
%\begin{equation}
%v\inww{0}{\pf+\pw-\mf}{} + \cfu\inww{\pf+\pw-\mf}{\pw}{} + \cfd\inww{\pw}{\pf+\pw}{} - \cfd\pf.
%\end{equation}
%\end{proposition}

\proof{Proof.} For given values of $\pf$ and $\pw$, the optimal re-dispatch actions $\pfu, \pfd, \iw, \sw$ are determined by solving  the following optimization problem.
\begin{subequations}
\begin{align}
& \of{\pfu, \pfd, \iw, \sw}{\inw{\om}{}{\pp{v \sw + \cfu \pfu - \cfd \pfd}}} \\
& \acc{\sw + \pfu - \pfd + \iw = 0, \quad \aom} \\
& \ac{0 \leq \pfd \leq \pf, \quad \aom} \\
& \ac{0 \leq \pfu \leq \mf - \pf, \quad \aom} \\
& \ac{0 \leq \pw + \iw \leq \ww, \quad \aom} \\
& \ac{0 \leq \sw \leq \pf + \pw, \quad \aom}
\end{align}
\end{subequations}

Depending on whether the realized stochastic generation $\ww$ is higher or lower than the dispatched quantity $\pw$, the re-dispatch rule for each scenario $\omega$ is given by:
\begin{equation}
\text{ If } \ww \leq \pw
\begin{cases}
\pfu = \miin{\mf-\pf,\pw-\ww} \\
\pfd = 0 \\
\sw=\maax{0,\pw-\ww-\mf+\pf} \\
\end{cases}
\end{equation}
\begin{equation}
\text{ If } \ww > \pw
\begin{cases}
\pfu = 0 \\
\pfd = \miin{\pf,\ww-\pw} \\
\sw=0 \\
\end{cases}
\end{equation}

Therefore, the second-stage expected cost can be computed as:
\begin{align}
& \inw{\om}{}{\pp{v \sw + \cfu \pfu - \cfd \pfd}} = \nonumber \\
& = \inw{0}{\pw}{\pp{v \sw + \cfu \pfu - \cfd \pfd}} + \inw{\pw}{\infty}{\pp{v \sw + \cfu \pfu - \cfd \pfd}} = \nonumber \\
& = \inw{0}{\pw}{v\sw} + \inw{0}{\pw}{\cfu \pfu} - \inw{\pw}{\infty}{\cfd \pfd} = \nonumber \\
& = v \inw{0}{\pf+\pw-\mf}{\pp{\pw-\ww-\mf+\pf}} + \cfu\inw{0}{\pf+\pw-\mf}{\pp{\mf-\pf}} + \nonumber \\
& \quad + \cfu\inw{\pf+\pw-\mf}{\pw}{\pp{\pw-\ww}} - \cfd \inw{\pw}{\pf+\pw}{\pp{\ww-\pw}} - \cfd \inw{\pf+\pw}{\infty}{\pf} = \nonumber \\
& = v\pp{\pf+\pw-\mf}\F{\pf+\pw-\mf} - v \inw{0}{\pf+\pw-\mf}{\ww} + \nonumber \\
& \quad + \cfu \pp{\mf-\pf-\pw}\F{\pf+\pw-\mf} + \cfu\pw\F{\pw} -  \cfu \inw{\pf+\pw-\mf}{\pw}{\ww} + \nonumber \\
& \quad + \cfd\pp{\pf+\pw}\F{\pf+\pw} - \cfd\pw\F{\pw} - \cfd\pf - \cfd\inw{\pw}{\pf+\pw}{\ww} = \nonumber \\
& = v\inww{0}{\pf+\pw-\mf}{} + \cfu\inww{\pf+\pw-\mf}{\pw}{} + \cfd\inww{\pw}{\pf+\pw}{} - \cfd\pf \label{2stage}
\end{align}

where, in the last equality, we have used the integration-by-parts theorem, according to which $\int_{x_1}^{x_2}sf(s)ds = x_2\F{x_2} - x_1\F{x_1} - \int_{x_1}^{x_2}\F{s}ds$. \hfill \Halmos \endproof

\section{Proof of Theorem \ref{Prop:ESCF}} \label{Proof:ESCF}

The proof of Theorem \ref{Prop:ESCF} relies on the proposition presented and proven below.

\begin{proposition} \label{Prop:WindStoc}
Consider a power system as described in Definition \ref{PowerSystem} in which the dispatch of the inflexible unit $\pii$ is given. Let $\pt$ then denote the net load to be satisfied by the flexible and stochastic power generating units. The dispatch rule for the flexible and stochastic generation is given by:

\vspace{4mm}
\begin{center}
{\footnotesize
\begin{tabular}{|l|ll|l|}
  \hline
  Rule \# & $p_{W}$ & $p_{F}$ & applies if \\
  \hline
  1. & $\pt$ & 0 & $0 \leq \pt \leq l_1$ \\
  2. & $l_1$ & $\pt-l_1$ & $l_1 < \pt < l_1 + \mf$ \\
  3. & $\pt-\mf$ & $\mf$ & $l_1 + \mf \leq \pt$ \\
  \hline
\end{tabular}}
\end{center}
\vspace{4mm}

Likewise, the marginal production cost of the flexible-stochastic generation portfolio, denoted by $\ct$, writes as:
\begin{equation}
\ct =
\begin{cases}
\pp{v-\cfu}\F{\pt - \mf} + \cfu \F{\pt}  & \ifq 0 \leq \pt \leq l_1 \\
\cf + \pp{v-\cfu}\F{\pt-\mf} - \cfd\FF{\pt}  & \ifq l_1 < \pt < l_1 + \mf \\
\pp{v-\cfd}\F{\pt-\mf}+\cfd\F{\pt}  & \ifq l_1 + \mf \leq \pt  \\
\end{cases} \label{mc_flex_wind}
\end{equation}
\end{proposition}

\proof{Proof.} Using Proposition \ref{Prop:Cost2Stage}, the optimal forward dispatch of a portfolio of flexible and stochastic power generation is given as the solution to the following optimization problem:
\begin{subequations} \label{windflex}
\begin{align}
& \of{}{\zt = \cf\pf + v\inww{0}{\pt-\mf}{} + \cfu\inww{\pt-\mf}{\pt-\pf}{} + \cfd\inww{\pt-\pf}{\pt}{} - \cfd\pf} \label{windflex_of}\\
& \acc{0 \leq \pf \leq \mf: \pp{\underline{\gamma},\overline{\gamma}}} \label{windflex_c1}
\end{align}
\end{subequations}
\noindent where $\underline{\gamma},\overline{\gamma}$ are dual variables and $\pw = \pt - \pf$. Problem \eqref{windflex} is a convex optimization problem that satisfies a Slater condition and therefore, the KKT conditions below are necessary and sufficient for optimality.
\begin{subequations} \label{kkt_wf}
\begin{align}
& \der{\zt}{\pf} - \underline{\gamma} + \overline{\gamma} = 0 \label{kkt_wf_1}\\
&0 \leq  \pf \perp \underline{\gamma} \geq 0 \label{kkt_wf_2}\\
& 0 \leq \pp{\mf - \pf} \perp \overline{\gamma} \geq 0 \label{kkt_wf_3}
\end{align}
\end{subequations}
\noindent where
\begin{equation}
\der{\zt}{\pf} = \pp{\cf-\cfd} - \pp{\cfu-\cfd}\F{\pt-\pf}
\end{equation}

We now determine the optimal solution to problem \eqref{windflex} by exhaustively enumerating the points that satisfy the optimality conditions \eqref{kkt_wf} as follows:

a) $\pf=0$
\begin{align*}
&\eqref{kkt_wf_3} \rightarrow \overline{\gamma}=0 \\
& \eqref{kkt_wf_1} \rightarrow \der{\zt}{\pf} \geq 0 \implies \F{\pt} \leq \nvs
\end{align*}

b) $0<\pf<\mf$
\begin{align*}
& \eqref{kkt_wf_2} \rightarrow \underline{\gamma} = 0 \\
& \eqref{kkt_wf_3} \rightarrow \overline{\gamma} = 0 \\
& \eqref{kkt_wf_1} \rightarrow \der{\zt}{\pf} = 0 \implies \F{\pt-\pf} = \nvs
\end{align*}

c) $\pf=\mf$
\begin{align*}
& \eqref{kkt_wf_2} \rightarrow \underline{\gamma} = 0 \\
& \eqref{kkt_wf_1} \rightarrow \F{\pt-\mf} \geq \nvs
\end{align*}

Since $\pw = \pt - \pf$, the solutions of cases a)--c) can be summarized as the following dispatch rule:
\begin{equation}
\pw =
\begin{cases}
\pt  & \ifq 0 \leq \pt \leq l_1 \\
l_1  & \ifq l_1 < \pt < l_1 + \mf \\
\pt - \mf &  \ifq l_1 + \mf \leq \pt  \\
\end{cases} \label{rule1}
\end{equation}

\noindent where $l_1=\Fi{\nvs}$. The marginal production cost of the flexible-stochastic generation portfolio $\ct$ is equal to $\der{\zt}{\pt}$, that is,
\begin{equation}
\ct =
\begin{cases}
\pp{v-\cfu}\F{\pt - \mf} + \cfu \F{\pt}  & \ifq 0 \leq \pt \leq l_1 \\
\cf + \pp{v-\cfu}\F{\pt-\mf} - \cfd\FF{\pt}  & \ifq l_1 < \pt < l_1 + \mf \\
\pp{v-\cfd}\F{\pt-\mf}+\cfd\F{\pt}  & \ifq l_1 + \mf \leq \pt  \\
\end{cases} \label{mc_flex_wind}
\end{equation}

Note that the functions $\pw\pp{\pt}$ and $\ct\pp{\pt}$ are increasing and continuous on $\pt$.
\hfill \Halmos \endproof

\begin{remark}[News-vendor solution] The characteristic constant $l_1$ in~\eqref{rule1} and~\eqref{mc_flex_wind} can be interpreted as the solution to the classical news-vendor problem~\citep{Raiffa}:
\begin{equation}
q_{opt}=\Fi{\dfrac{p-c}{p+h}}
\end{equation} \label{newsvendor}
\noindent where $q_{opt}$ is the optimal stocking quantity of the news-vendor, $\F{\cdot}$ is the cumulative distribution function of the demand to be satisfied, $c$ is the variable production cost, and $p$ and $h$ correspond to the penalty cost of unsatisfied orders and the inventory holding cost, respectively. The analogy works, thus, as follows: The variable cost of the stochastic power production is zero, i.e., $c = 0$; $p = \cf-\cfd$ represents the penalty cost of dispatching less stochastic power capacity than its eventual power production, since the consequent power surplus is to be compensated for by a decrease in the flexible power generation (downward regulation); finally, $h = \cfu-\cf$ provides the marginal cost of dispatching more stochastic power capacity than its eventual real-time power production, because the consequent generation deficit is to be covered with an increase in the flexible power generation (upward regulation). Therefore,

\begin{equation}
q_{opt}=\Fi{\dfrac{p-c}{p+h}} =  \Fi{\frac{\cf - \cfd}{\cf - \cfd + \cfu  -\cf}}= \Fi{\nvs}= l_1.
\end{equation} \label{newsvendor}
\end{remark}

\begin{remark}[Discrete probability distribution]
Proposition \ref{Prop:WindStoc} assumes a continuous probability distribution for the uncertain electricity supply, which implies that objective function \eqref{windflex_of} is differentiable. If the uncertain supply is, in contrast, characterized by a discrete probability distribution, the cumulative distribution function  $\F{\cdot}$ is stepwise and thus, objective function \eqref{windflex_of} becomes nondifferentiable. In that case, the sub-derivative of the objective function should be used instead to formulate the KKT optimality conditions as follows:
\begin{subequations} \label{kkt_wf2}
\begin{align}
& \der{\zt^-}{\pf} + \overline{\gamma} - \underline{\gamma} \leq 0 \leq  \der{\zt^+}{\pf} + \overline{\gamma} - \underline{\gamma} \label{kkt_wf2_1}\\
& 0 \leq \pf \perp \underline{\gamma} \geq 0 \label{kkt_wf2_2}\\
& 0 \leq \pp{\mf - \pf} \perp \overline{\gamma} \geq 0 \label{kkt_wf2_3}
\end{align}
\end{subequations}
\noindent where
\begin{subequations}
\begin{align}
& \der{\zt^+}{\pf} = \cf - \cfd - \pp{\cfu-\cfd}\F{\pt-\pf} \\
& \der{\zt^-}{\pf} = \cf - \cfd - \pp{\cfu-\cfd}\pp{\F{\pt-\pf}-\f{\pt-\pf}}
\end{align}
\end{subequations}
Note that $f(\cdot)$ should be interpreted here as the probability mass function.
We analyze next the different points satisfying the optimality conditions \eqref{kkt_wf2}:

a) $\pf=0$
\begin{align*}
&\eqref{kkt_wf2_3} \rightarrow \overline{\gamma} = 0\\
&\eqref{kkt_wf2_1} \rightarrow \der{\zt^+}{\pf}  - \underline{\gamma} \geq 0 \implies \der{\zt^+}{\pf} \geq 0 \implies \F{\pt} \leq \nvs
\end{align*}

b) $0<\pf<\mf$
\begin{align*}
&\eqref{kkt_wf2_2} \rightarrow \underline{\gamma} = 0\\
&\eqref{kkt_wf2_3} \rightarrow \overline{\gamma} = 0\\
& \eqref{kkt_wf2_1} \rightarrow \nvs \leq \F{\pt-\pf} \leq \nvs + \f{\pt-\pf}
\end{align*}

c) $\pf=\mf$
\begin{align*}
&\eqref{kkt_wf2_2} \rightarrow \underline{\gamma} = 0\\
& \eqref{kkt_wf2_1} \rightarrow  \der{\zt^-}{\pf} + \overline{\gamma} \leq 0 \implies \cf - \cfd - \pp{\cfu-\cfd}\pp{\F{\pt-\pf} - \f{\pt-\pf}} + \overline{\gamma} \leq 0 \implies \\
& \qquad \implies \F{\pt-\mf} \geq \nvs
\end{align*}

Since $\Fi{\cdot}$ is the generalized inverse distribution function, and $l_1$ is defined as $\Fi{\nvs}$, the dispatch rule for a discrete probability distribution of the uncertain supply coincides with \eqref{rule1}.

\end{remark}

Using Proposition \ref{Prop:WindStoc}, the proof of Theorem \ref{Prop:ESCF} proceeds as follows:

\proof{Proof.} This proof deals with the optimal dispatch of an inflexible power unit with capacity $\mi>0$ and marginal cost $\ci$ and the flexible-stochastic generation portfolio of Proposition \ref{Prop:WindStoc}. For a total system load denoted by $l$, this optimal dispatch is determined by solving the following optimization problem:
\begin{subequations} \label{inflexwindflex}
\begin{align}
& \of{\pii \geq 0,\pt \geq 0}{z = \ci\pii + \int_{0}^{\pt} \ct\pp{x}dx} \\
& \pii + \pt = l: \tau \\
& \pii \leq \mi: \phi
\end{align}
\end{subequations}

Since the integral of an increasing function is a convex function, problem \eqref{inflexwindflex} is a convex optimization problem and therefore, the KKT conditions below are necessary and sufficient for optimality.
\begin{subequations} \label{kkt_iwf}
\begin{align}
& 0 \leq  \pp{\ci - \tau + \phi} \perp \pii \geq 0 \label{kkt_iwf_1}\\
& 0 \leq \pp{\ct(\pt) - \tau} \perp \pt \geq 0 \label{kkt_iwf_2}\\
& 0 \leq \pp{\mi-\pii} \perp \phi \geq 0 \label{kkt_iwf_3}\\
& \pii + \pt = l \label{kkt_iwf_4}
\end{align}
\end{subequations}

The solution to problem \eqref{inflexwindflex} is then obtained by exhaustively examining the points that satisfy the optimality conditions \eqref{kkt_iwf} as follows:

a) $\pii = 0\andq \pt = 0$
\begin{align*}
& \eqref{kkt_iwf_4} \rightarrow \text{only feasible if} \; l=0 \\
& \eqref{kkt_iwf_3} \rightarrow \phi = 0 \\
& \eqref{kkt_iwf_1} \rightarrow \tau \leq \ci \\
& \eqref{kkt_iwf_2} \rightarrow \tau \leq \ct(0)
\end{align*}

b) $\pii = 0\andq \pt >0 $
\begin{align*}
& \eqref{kkt_iwf_4} \rightarrow \pt = l \\
&\eqref{kkt_iwf_3} \rightarrow \phi = 0 \\
& \begin{rcases}
        \eqref{kkt_iwf_1} \rightarrow \tau = \ct(\pt) \\
        \eqref{kkt_iwf_2} \rightarrow \tau \leq \ci \\
        \end{rcases} \ct(\pt) \leq \ci
\end{align*}

c) $\pii = \mi \andq \pt = 0$
\begin{align*}
&\eqref{kkt_iwf_4} \rightarrow \mi = l \\
&\eqref{kkt_iwf_3} \rightarrow \phi \geq 0 \\
&\begin{rcases}
        \eqref{kkt_iwf_1} \rightarrow \tau \geq \ci \\
        \eqref{kkt_iwf_2} \rightarrow \tau \leq \ct(0) \\
        \end{rcases} \ci \leq \tau \leq \ct(0)
\end{align*}

d) $\pii = \mi \andq \pt > 0 $
\begin{align*}
& \eqref{kkt_iwf_4} \rightarrow \pt = l - \mi\\
&\eqref{kkt_iwf_3} \rightarrow \phi \geq 0\\
&\begin{rcases}
        \eqref{kkt_iwf_1} \rightarrow \tau \geq \ci \\
        \eqref{kkt_iwf_2} \rightarrow \tau = \ct(\pt)
        \end{rcases} \ci \leq \ct(\pt)
\end{align*}

e) $ 0 < \pii < \mi \andq \pt = 0$
\begin{align*}
&\eqref{kkt_iwf_4} \rightarrow \pii = l \\
&\eqref{kkt_iwf_3} \rightarrow \phi = 0 \\
&\begin{rcases}
        \eqref{kkt_iwf_1} \rightarrow \tau = \ci \\
        \eqref{kkt_iwf_2} \rightarrow \tau \leq \ct(0) \\
        \end{rcases} \ci \leq \ct(0)
\end{align*}

f) $ 0 < \pii < \mi \andq \pt > 0 $
\begin{align*}
& \eqref{kkt_iwf_4} \rightarrow \pii + \pt = l \\
& \eqref{kkt_iwf_3} \rightarrow \phi = 0 \\
& \begin{rcases}
        \eqref{kkt_iwf_1} \rightarrow \tau = \ci \\
        \eqref{kkt_iwf_2} \rightarrow \tau = \ct(\pt) \\
        \end{rcases} \ci = \ct(\pt)
\end{align*}

Let $\ph$ denote the value of $\pt$ such that $\ct\pp{\ph} = \ci$. Note that the function $\ct(\pt)$ is increasing and continuous on $\pt$ and therefore, $\ph$ exists provided that $\ct(0)\leq\ci\leq\ct(\infty)$. If $\ci<\ct(0)$, then we set $\ph=0$. Similarly, if $\ci>\ct(\infty)$, we assign the value $\infty$ to $\ph$. This way the solutions obtained in cases a)--f) above can be summarized in the following dispatch rule:
\begin{equation}
\pt =
\begin{cases}
l     &\ifq 0 \leq l \leq \ph \\
\ph   &\ifq \ph \leq l \leq \ph + \mi \\
l-\mi &\ifq \ph + \mi \leq l \\
\end{cases} \label{rule2}
\end{equation}

\noindent where the dispatch of the inflexible power capacity is $\pii = l - \pt$. Note that the function $\ct(\pt)$ is piecewise and therefore,  $\ph$ can take on the following values
\begin{equation}
\ph =
\begin{cases}
0 \ifq \ci\leq \ct(0) \\
l_2 \ifq \ct(0) \leq \ci \leq \ct(l_1) \\
l_3 \ifq \ct(l_1) \leq \ci \leq \ct(l_1+\mf) \\
l_4 \ifq \ct(l_1+\mf) \leq \ci \leq \ct(\infty) \\
\infty \ifq \ct(\infty) \leq \ci \\
\end{cases} \label{phat}
\end{equation}

where
\begin{align*}
& l_2 := l : \pp{v-\cfu}\F{l - \mf} + \cfu \F{l} = \ci \\
& l_3 := l : \pp{v-\cfu}\F{l-\mf} - \cfd\FF{l} = \ci - \cf \\
& l_4 := l : \pp{v-\cfd}\F{l-\mf}+\cfd\F{l} = \ci
\end{align*}

Merging \eqref{rule1}, \eqref{rule2} and \eqref{phat} we obtain the following cases:

\vspace{4mm} \begin{center}
{\footnotesize
\begin{tabular}{|l|lll|l|l|}
\hline
Rule \# & $\pw$ & $\pii$ & $\pf$ & \multicolumn{2}{c|}{applies if} \\
\hline
1. & $0$ & $l$ & $0$ & \multirow{4}{*}{$\ci\leq \ct(0)$} & $0 \leq l \leq \mi$ \\
2. & $l-\mi$ & $\mi$ & $0$ &  & $\mi \leq l \leq l_1 + \mi$ \\
3. & $l_1$ & $\mi$ & $l-\l_1-\mi$ &  & $l_1 + \mi \leq l \leq l_1 + \mi + \mf$ \\
4. & $l-\mi-\mf$ & $\mi$ & $\mf$ &  & $l_1 + \mi + \mf \leq l$ \\
\hline
5. & $l$ & $0$ & $0$ & \multirow{5}{*}{$\ct(0) \leq \ci \leq \ct(l_1)$} & $0\leq l \leq l_2$ \\
6. & $l_2$ & $l-l_2$ & $0$ & & $l_2\leq l \leq l_2+\mi$ \\
7. & $l-\mi$ & $\mi$ & $0$ & & $l_2+\mi\leq l \leq l_1+\mi$ \\
8. & $l_1$ & $\mi$ & $l-l_1-\mi$ & & $l_1+\mi\leq l \leq l_1+\mi+\mf$ \\
9. & $l-\mi-\mf$ & $\mi$ & $\mf$ & & $l_1+\mi+\mf\leq l$ \\
\hline
10. & $l$ & $0$ & $0$ & \multirow{5}{*}{$\ct(l_1) \leq \ci \leq \ct(l_1+\mf)$} & $0\leq l \leq l_1$ \\
11. & $l_1$ & $0$ & $l-l_1$ & & $l_1\leq l \leq l_3$ \\
12. & $l_1$ & $l-l_3$ & $l_3-l_1$ & & $l_3\leq l \leq l_3+\mi$ \\
13. & $l_1$ & $\mi$ & $l-l_1-\mi$ & & $l_3+\mi\leq l \leq l_1+\mf+\mi$ \\
14. & $l-\mi-\mf$ & $\mi$ & $\mf$ & & $l_1+\mf+\mi\leq l$ \\
\hline
15. & $l$ & $0$ & $0$ & \multirow{5}{*}{$\ct(l_1+\mf) \leq \ci \leq \ct(\infty)$} & $0 \leq l \leq l_1$ \\
16. & $l_1$ & $0$ & $l-l_1$ & & $l_1 \leq l \leq l_1+\mf$\\
17. & $l-\mf$ & $0$ & $\mf$ & & $l_1+\mf \leq l \leq l_4$\\
18. & $l_4-\mf$ & $l-l_4$ & $\mf$ & & $l_4 \leq l \leq l_4+\mi$\\
19. & $l-\mi-\mf$ & $\mi$ & $\mf$ & & $l_4+\mi \leq l $\\
\hline
20. & $l$ & $0$ & $0$ & \multirow{3}{*}{$\ct(\infty) \leq \ci$} & $0 \leq l \leq l_1$\\
21. & $l_1$ & $0$ & $l-l_1$ & & $l_1 \leq l \leq l_1+\mf$\\
22. & $l-\mf$ & $0$ & $\mf$ & & $l_1+\mf \leq l$\\
\hline
\end{tabular}} \end{center}
\vspace{4mm}

Next the different dispatch rules above are recast as a function of the system characteristics constants $l_1,l_2,l_3,l_4$ only. To do so, we note that if $l_2,l_3,l_4$ do not exit, their values are assigned to infinity.
%Also, we need to impose that $l_3 \geq l_1$ and $l_4 \geq \l_1 + \mf$.

Rules 1-4 only apply if the marginal cost of the stochastic-flexible portfolio for $\pt=0$ is higher than the marginal cost of the inflexible power unit $\ci$. Note also that the condition $\ci \leq \ct(l_1)$ is equivalent to $l_1 \geq l_2$ and thus we can jointly reformulate rules 1-9 as follows:

\vspace{4mm} \begin{center}
{\footnotesize
\begin{tabular}{|lll|l|l|l|l|}
  \hline
 $p_{W}$ & $p_{I}$ & $p_{F}$ & \multicolumn{4}{c|}{applies if}   \\
  \hline
 $l$ & $0$ & $0$                                                             & \multirow{5}{1.2cm}{$l_1 \geq l_2$}  & \multicolumn{3}{l|}{$0 \leq l \leq l_2$} \\
 $l_2$ &$l - l_2$ & $0$                                                       &                                      & \multicolumn{3}{l|}{$l_2 < l \leq \mi + l_2$}\\
 $l-\mi$ & $\mi$ & $0$                                 &                                      & \multicolumn{3}{l|}{$\mi + l_2 < l \leq \mi + l_1$}\\
 $l_1$ & $\mi$ &$l - l_1 -\mi$                        &                                      & \multicolumn{3}{l|}{$\mi + l_1 < l \leq \mf + \mi + l_1$}\\
 $l - \mf - \mi$ &$\mi$& $\mf$ &                                      & \multicolumn{3}{l|}{$l > \mf + \mi + l_1$}\\
\hline
\end{tabular}} \end{center}
\vspace{4mm}

Rules 10,15 and 20 can be easily merged as

\vspace{4mm} \begin{center}
{\footnotesize
\begin{tabular}{|lll|l|l|l|l|}
  \hline
  $p_{W}$ & $p_{I}$ & $p_{F}$ & \multicolumn{4}{c|}{applies if}   \\
  \hline
  $l$ & $0$ & $0$                                                             & $l_1 < l_2$ & \multicolumn{3}{l|}{$0 \leq l \leq l_1$}                                  \\
 \hline
\end{tabular}} \end{center}
\vspace{4mm}

The condition $\ct(l_1) \leq \ci \leq \ct(l_1+\mf)$ can be equivalently formulated as $l_1 < l_2$ and $l_3 \leq l_1 + \mf$, which allows expressing rules 11-14 as:

\vspace{4mm} \begin{center}
{\footnotesize
\begin{tabular}{|lll|l|l|l|}
  \hline
  $p_{W}$ & $p_{I}$ & $p_{F}$ & \multicolumn{3}{c|}{applies if}   \\
  \hline
 $l_1$ &$0$ & $l-l_1$  & \multirow{4}{1.2cm}{$l_1 < l_2$}  & \multirow{4}{2.5cm}{$l_3 \leq l_1 + \mf$} & $l_1 \leq l \leq l_3$ \\
 $l_1$ &$l - l_3$ & $l_3 - l_1$ &  & & $l_3 < l \leq \mi + l_3$ \\
 $l_1$ & $\mi$ & $l - l_1 - \mi$ & & & $l_3 + \mi < l \leq \mi + \mf + l_1$ \\
$l - \mi -\mf$ & $\mi$ & $\mf$ & & & $l > \mi + \mf + l_1$ \\
 \hline
\end{tabular}} \end{center}
\vspace{4mm}

Rules 16 and 21 can also be merged as:

\vspace{4mm} \begin{center}
{\footnotesize
\begin{tabular}{|lll|l|l|l|l|}
  \hline
  $p_{W}$ & $p_{I}$ & $p_{F}$ & \multicolumn{4}{c|}{applies if}   \\
  \hline
  $l_1$ & $0$ & $l-l_1$                                                       &              \multirow{1}{1.2cm}{$l_1 < l_2$}                        & \multirow{1}{3cm}{$l_3 > l_1 + \mf$}      & \multicolumn{2}{l|}{$l_1 \leq l \leq \mf + l_1$} \\
  \hline
\end{tabular}} \end{center}
\vspace{4mm}

Likewise, rules 17-19 and 22 can be rewritten as:

\vspace{4mm} \begin{center}
{\footnotesize
\begin{tabular}{|lll|l|l|l|l|}
  \hline
   $p_{W}$ & $p_{I}$ & $p_{F}$ & \multicolumn{4}{c|}{applies if}   \\
  \hline
 $l - \mf$ & $0$ & $\mf$                               &                                     \multirow{4}{1.2cm}{$l_1 < l_2$} &                \multirow{4}{3cm}{$l_3 > l_1 + \mf$}                                      & \multirow{3}{1.2cm}{$\exists l_4$ $(c_I \leq v)$}   &  $\mf + l_1 \leq l \leq l_4$\\
 $l_4 - \mf$ & $l - l_4$ & $\mf$                       &                                      &                                                      &                       & $l_4 \leq l \leq l_4 + \mi$\\
 $l - \mf - \mi$ & $\mi$ & $\mf$ &                                      &                                                      &                      & $l > l_4 + \mi$\\
\cline{6-7}
 $l - \mf$ & $0$ & $\mf$                               &                                      &                                                      &  \multirow{1}{1.2cm}{$\nexists l_4$}      &  $l > \mf + l_1$\\
  \hline
\end{tabular}} \end{center}
\vspace{4mm}

The final dispatch rule can be thus summarized as:

\vspace{4mm} \begin{center}
{\footnotesize
\begin{tabular}{|l|lll|l|l|l|l|}
  \hline
  Rule \# & $p_{W}$ & $p_{I}$ & $p_{F}$ & \multicolumn{4}{c|}{applies if}   \\
  \hline
 1. & $l$ & $0$ & $0$                                                             & \multirow{5}{1.2cm}{$l_1 \geq l_2$}  & \multicolumn{3}{l|}{$0 \leq l \leq l_2$} \\
 2. & $l_2$ &$l - l_2$ & $0$                                                       &                                      & \multicolumn{3}{l|}{$l_2 < l \leq \mi + l_2$}\\
 3. & $l-\mi$ & $\mi$ & $0$                                 &                                      & \multicolumn{3}{l|}{$\mi + l_2 < l \leq \mi + l_1$}\\
 4. & $l_1$ & $\mi$ &$l - l_1 -\mi$                        &                                      & \multicolumn{3}{l|}{$\mi + l_1 < l \leq \mf + \mi + l_1$}\\
 5. & $l - \mf - \mi$ &$\mi$& $\mf$ &                                      & \multicolumn{3}{l|}{$l > \mf + \mi + l_1$}\\
\hline
 6. & $l$ & $0$ & $0$                                                             & \multirow{9}{1.2cm}{$l_1 < l_2$}     & \multicolumn{3}{l|}{$0 \leq l \leq l_1$}                                  \\
 \cline{6-8}
 7. & $l_1$ &$0$ & $l-l_1$                                                       &                                      & \multirow{4}{3cm}{$l_3 \leq l_1 + \mf$} & \multicolumn{2}{l|}{$l_1 \leq l \leq l_3$} \\
 8. & $l_1$ &$l - l_3$ & $l_3 - l_1$                                             &                                      &                                                      & \multicolumn{2}{l|}{$l_3 < l \leq \mi + l_3$} \\
 9. & $l_1$ & $\mi$ & $l - l_1 - \mi$                       &                                      &                                                      & \multicolumn{2}{l|}{$l_3 + \mi < l \leq \mi + \mf + l_1$} \\
10. & $l - \mi -\mf$ & $\mi$ & $\mf$  &                                      &                                                      & \multicolumn{2}{l|}{$l > \mi + \mf + l_1$} \\
\cline{1-4}\cline{6-8}
11. & $l_1$ & $0$ & $l-l_1$                                                       &                                      & \multirow{5}{3cm}{$l_3 > l_1 + \mf$}      & \multicolumn{2}{l|}{$l_1 \leq l \leq \mf + l_1$} \\
\cline{7-8}
12. & $l - \mf$ & $0$ & $\mf$                               &                                      &                                                      & \multirow{3}{1.2cm}{$\exists l_4$ $(c_I \leq v)$}   &  $\mf + l_1 \leq l \leq l_4$\\
13. & $l_4 - \mf$ & $l - l_4$ & $\mf$                       &                                      &                                                      &                       & $l_4 \leq l \leq l_4 + \mi$\\
14. & $l - \mf - \mi$ & $\mi$ & $\mf$ &                                      &                                                      &                      & $l > l_4 + \mi$\\
\cline{7-8}
15. & $l - \mf$ & $0$ & $\mf$                               &                                      &                                                      &  \multirow{1}{1.2cm}{$\nexists l_4$}      &  $l > \mf + l_1$\\
  \hline
\end{tabular}} \end{center}
\vspace{4mm}

\hfill \Halmos \endproof

\begin{remark}[Discrete probability distribution]
The same sub-differential analysis used in Proposition \ref{Prop:WindStoc} can be applied here so that the dispatch rule above is also valid for the case in which the uncertain power supply is modeled by a discrete probability distribution. In such a case, it suffices to redefine the constants $l_2,l_3,l_4$ as:
\begin{align*}
l_2 &:= \underset{l \geq 0}{\rm min} \enskip l: \left(v-c_F^+\right)F\left(l-\mf\right) + c_F^+F(l) \geq c_I \\
l_3 &:= \underset{l \geq l_1}{\rm min} \enskip l: c_F + \left(v-c_F^+\right)F\left(l-\mf\right) - c_F^- \left(1-F(l)\right) \geq c_I \\
l_4 &:= \underset{l \geq l_1 + \mf}{\rm min} \enskip l: \left(v-c_F^-\right)F\left(l-\mf\right) + c_F^- F(l) \geq c_I
\end{align*}
\end{remark}

\section{Proof of Corollary \ref{Corr:CapacityAdequate}} \label{Proof:CapacityAdequate}

\proof{Proof.} %
First, we show that, under the conditions stated in this corollary, $l_3 \leq l_1 + \mf$, and therefore, rules 11--15 of the stochastic dispatch solution do \emph{not} apply. For this purpose, consider expression~\eqref{L3}, which defines constant $l_3$, and note that $c_F + \left(v-c_F^+\right)F\left(\mf-\mf\right) - c_F^- \left(1-F(\mf)\right) = \cf + (v-\cfu) \F{0}$, because $\F{\mf} = 1$ given that $\mf \geq \mw$. Furthermore, it holds that  $\cf + (v-\cfu) \F{0} > \ci$, since $\cf > \ci$ and $(v-\cfu) \F{0} \geq 0$. Consequently, we have that $l_3 \leq \mf \leq l_1 + \mf$.

Now we prove that $l_3 = max\left(l_1, A\right)$ with $A = F^{-1}\left(1-\frac{c_{F}-c_{I}}{c_{F}^{-}}\right)$ is optimal for the minimization problem~\eqref{L3}. By construction, $l_3 \leq \mw$. Note that
$c_F + \left(v-c_F^+\right)F\left(A-\mf\right) - c_F^- \left(1-F(A)\right) - \ci = \left(v-c_F^+\right)F\left(A-\mf\right) \geq 0$ and that $\F{x} = 0$ for all $x < 0$. Hence, either $l_3 = A$ if $A \geq l_1$ or $l_3 = l_1$ otherwise.

We can proceed analogously to show that $l_2 = F^{-1}\left(\frac{c_{I}}{c_{F}^{+}}\right)$. Let denote constant $F^{-1}\left(\frac{c_{I}}{c_{F}^{+}}\right)$ by $B$. By construction $0 \leq B \leq \mw$. Consider minimization problem~\eqref{L2}. It holds that $\left(v-\cfu\right) \F{B-\mf} + \cfu\F{B}-\ci = \left(v-\cfu\right) \F{B-\mf} \geq 0$. Since $B-\mf \leq 0$ and $\F{x} = 0$ for all $x < 0$, $l_2 = B$ is optimal for \eqref{L2}.

Finally, the result corresponding to $l\leq\mi$ trivially follows by isolating those ranges of system load for which $l$ must be lower than $\mi$. \hfill \Halmos \endproof

\section{Proof of Theorem \ref{Prop:VB}} \label{Proof:VB}

\proof{Proof.} Let us start identifying the optimal solution to the forward dispatch model \eqref{Clearing_DA} for a given value of the virtual bid volume $\pv$. The merit-order dispatch of the inflexible, flexible and stochastic power capacity is provided in the table below. Note that the case $l > \uw + \pv + \mi + \mf$ makes problem \eqref{Clearing_DA} infeasible.

\vspace{4mm}
\begin{center}
{\footnotesize
\begin{tabular}{|c|llll|l|}
\hline
case & $\pf$ & $\pii$ & $\pw$ & $\lf$ & applies if \\
\hline
a) & $0$ & $0$ & $l-\pv$ & $0$ & $l < \uw + \pv$ \\
b) & $0$ & $0$ & $\uw$ & $[0,\ci]$ & $l = \uw + \pv$ \\
c) & $0$ & $l - \uw - \pv$ & $\uw$ & $\ci$ & $\uw+\pv < l < \uw+\pv+\mi$ \\
d) & $0$ & $\mi$ & $\uw$ & $[\ci,\cf]$ & $l = \uw + \pv + \mi$ \\
e) & $l - \uw - \pv - \mi$ & $\mi$ & $\uw$ & $\cf$ & $\uw + \pv + \mi < l < \uw + \pv + \mi + \mf$ \\
f) & $\mf$ & $\mi$ & $\uw$ & $[\cf,\infty]$ & $l = \uw + \pv + \mi + \mf$ \\
\hline
\end{tabular}} \end{center}
\vspace{4mm}

Observe that in cases a), b) or c), the flexible generation technology is not dispatched and therefore, the balancing market-clearing problem \eqref{Clearing_BM} simplifies to:
\begin{subequations}
\begin{align}
& \of{\pfu, \iw, \sw}{\cfu\pfu + v\sw} \\
& \acc{\pfu + \iw + \sw + \iv = 0: \lb} \\
& \ac{0 \leq \pfu \leq \mf} \\
& \ac{0 \leq \pw + \iw \leq \ww} \\
& \ac{0 \leq \sw \leq l}
\end{align}
\end{subequations}

whose solution boils down to a merit-order-based dispatch, that is,

\vspace{4mm}
{\footnotesize
\begin{tabular}{|lllll|l|}
\hline
$\pfu$ & $\pfd$ & $\iw$ & $\sw$ & $\lb$ & applies if \\
\hline
$0$ & $0$ & $\pv$ & $0$ & $0$ & $\pw + \pv < \ww$\\
$0$ & $0$ & $\pv$ & $0$ & $[0,\cfu]$ & $\ww = \pw + \pv$\\
$\pw + \pv - \ww$ & $0$ & $\ww - \pw$ & $0$ & $\cfu$ & $\pw + \pv - \mf < \ww < \pw + \pv$\\
$\pw + \pv - \ww$ & $0$ & $\ww - \pw$ & $0$ & $[\cfu,v]$ & $\ww = \pw + \pv - \mf$\\
$\mf$ & $0$ & $\ww - \pw$ & $\pw+\pv-\ww-\mf$ & $v$ & $\ww < \pw + \pv -  \mf$\\
\hline
\end{tabular}}\\

The expected balancing price for these cases is computed as:
\begin{align*}
& \inw{\om}{}{\lb} = \inw{0}{\pw+\pv-\mf}{v} + \inw{\pw + \pv - \mf}{\pw + \pv}{\cfu} = v\F{\pw+\pv-\mf}+ \nonumber \\
& \qquad + \cfu\pp{\F{\pw+\pv}-\F{\pw+\pv-\mf}} = \pp{v-\cfu}\F{\pw+\pv-\mf} + \cfu\F{\pw+\pv}
\end{align*}

In cases d), e) and f) of the forward dispatch, $\pw=\uw$ and therefore, the balancing market-clearing problem~\eqref{Clearing_BM} becomes:

\begin{subequations}
\begin{align}
& \of{\pfu,\pfd,\iw,\sw}{\cfu\pfu - \cfd\pfd + v\sw} \\
& \acc{\pfu - \pfd + \iw + \sw + \iv = 0: \lb} \\
& \ac{0 \leq \pfd \leq \pf} \\
& \ac{0 \leq \pfu \leq \mf - \pf} \\
& \ac{0 \leq \uw + \iw \leq \ww} \\
& \ac{0 \leq \sw \leq l}
\end{align}
\end{subequations}

whose solution is given by:

\vspace{4mm}
{\footnotesize
\begin{tabular}{|lllll|l|}
\hline
$\pfu$ & $\pfd$ & $\iw$ & $\sw$ & $\lb$ & applies if \\
\hline
$0$ & $l-\mi-\uw-\pv$ & $l - \mi - \uw$ & $0$ & $0$ & $l-\mi < \ww$\\
$0$ & $l - \mi - \uw - \pv$ & $l - \mi - \uw$ & $0$ & $[0,\cfd]$ & $ \ww = l-\mi$ \\
$0$ & $\ww-\uw-\pv$ & $\ww-\uw$ & $0$ & $\cfd$ & $\uw+\pv < \ww < l-\mi $\\
$0$ & $0$ & $\pv$ & $0$ & $[\cfd,\cfu]$ & $\ww = \uw + \pv $\\
$\uw + \pv - \ww$ & $0$ & $\ww - \uw$ & $0$ & $\cfu$ & $l-\mi-\mf < \ww < \uw + \pv $\\
$\uw + \pv - l + \mi + \mf$ & $0$ & $l - \mi - \mf - \uw$ & $0$ & $ [\cfu,v]$  & $\ww=l-\mi-\mf $ \\
$\mf-l+\uw+\pv+\mi$ & $0$ & $\ww-\uw$ & $l-\mf-\mi-\ww$ & $v$ & $ \ww < l-\mi-\mf $\\
\hline
\end{tabular}}\\

Hence, the expected balancing price in cases d), e) and f) is calculated as:
\begin{align*}
& \inw{\om}{}{\lb} = \inw{0}{l-\mi-\mf}{v} + \inw{l-\mi-\mf}{\uw+\pv}{\cfu} + \inw{\uw+\pv}{l-\mi}{\cfd} + \inw{l-\mi}{\infty}{0} = \nonumber \\
& \qquad = v\F{l-\mi-\mf} + \cfu\pp{\F{\uw+\pv}-\F{l-\mi-\mf}} + \cfd\pp{\F{l-\mi}-\F{\uw+\pv}} = \nonumber \\
& \qquad = \pp{v-\cfu}\F{l-\mi-\mf} + \pp{\cfu-\cfd}\F{\uw+\pv} + \cfd\F{l-\mi}
\end{align*}

Finally, the KKT conditions of the arbitrager's problem \eqref{Arbitrager} imply that:
\begin{equation}
\der{}{\pv} \pp{\lf\pv - \pv\inw{\om}{}{\lb}} = 0 \implies \lf = \inw{\om}{}{\lb}
\end{equation}

That is, the strategy of the arbitrager is to place a zero-price virtual bid $\pv$ in the forward market and to repurchase or resell the same amount in the balancing market so that the forward price equals the expected balancing price.

The solution to the short-run equilibrium problem must be a point that is simultaneously optimal for the arbitrager's problem and the clearing problems of the forward and the balancing markets. Therefore, in the short-run equilibrium solution the forward price $\lf$ must be equal to the expected balancing price $\inw{\om}{}{\lb}$, which we denote by $\lbh$. Below we analyze the solution to the short-run equilibrium problem for each of the cases a) to f) included in the first table of this proof.

a) \begin{equation*}
\begin{rcases}
\lf = 0 \\
\lbh = (v-\cfu)\F{l-\mf}+\cfu\F{l}
\end{rcases} l=0
\end{equation*}

b) \begin{equation*}
\begin{rcases}
    0 \leq \lf \leq \ci \\
    \lbh = \pp{v-\cfu}\F{l-\mf} + \cfu\F{l}
    \end{rcases} 0 \leq l \leq l_2
\end{equation*}

c) \begin{equation*}
\begin{rcases} \begin{rcases}
    \lf = \ci \\
    \lbh = \pp{v-\cfu}\F{\uw+\pv-\mf} + \cfu\F{\uw+\pv}
    \end{rcases} \uw + \pv = l_2 \\
\qquad \qquad \qquad \qquad \qquad \qquad \qquad \qquad \uw+\pv < l < \uw+\pv+\mi \end{rcases} l_2 < l < l_2+\mi
\end{equation*}

d) \begin{equation*}
\begin{rcases}
    \ci \leq \lf \leq \cf \\
    \lbh = \pp{v-\cfu}\F{l-\mi-\mf} + \cfu\F{l-\mi}
    \end{rcases} l_2 \leq l - \mi \leq l_5
\end{equation*}

e) \begin{equation*}
\begin{rcases} \begin{rcases}
    \lf = \cf \\
    \lbh = \pp{v-\cfu}\F{l-\mi-\mf} + \pp{\cfu-\cfd}\F{\uw+\pv} + \cfd\F{l-\mi}
    \end{rcases} \implies \\
\qquad \implies  \F{\uw+\pv}= \dfrac{\cf-\pp{v-\cfu}\F{l-\mi-\mf}-\cfd\F{l-\mi}}{\cfu-\cfd} \\
\uw + \pv + \mi < l \leq \uw + \pv + \mi + \mf \implies \\
\qquad \implies  l - \mi - \mf \leq \uw + \pv \leq l - \mi \implies \\
\qquad \implies \F{l-\mi-\mf} < \F{\uw+\pv} \leq \F{l- \mi} \end{rcases} l_5 < l-\mi < l_6
\end{equation*}

f) \begin{equation*}
\begin{rcases}
    \cf \leq \lf  \\
    \lbh = \pp{v-\cfd}\F{l-\mi-\mf} + \cfd\F{l-\mi}
    \end{rcases} l-\mi \geq l_6
\end{equation*}

\noindent where
%
%\begin{align*}
%& l_2 := \pt : \pp{v-\cfu}\F{\pt - \mf} + \cfu \F{\pt} = \ci \\
%& l_5 := \pt : \pp{v-\cfu}\F{\pt - \mi - \mf} + \cfu \F{\pt - \mi} = \cf \\
%& l_6 := \pt : \pp{v-\cfd}\F{\pt- \mi - \mf}+\cfd\F{\pt - \mi} = \cf \\
%& l_7(l) = \Fi{\dfrac{\cf - \pp{v-\cfu}\F{l-\mi-\mf} - \cfd\F{l-\mi}}{\cfu - \cfd}}\\
%\end{align*}
\begin{align*}
l_2 &:= \underset{l \geq 0}{\rm min} \enskip l: \left(v-c_F^+\right)F\left(l-\mf\right) + c_F^+F(l) \geq c_I\\
l_5 &:= \underset{l \geq 0}{\rm min} \enskip l: \left(v-c_F^+\right)F\left(l-\mf- \mi\right) + c_F^+F(l- \mi) \geq c_F\\
l_6 &:= \underset{l \geq 0}{\rm min} \enskip l: \left(v-c_F^-\right)F\left(l-\mf- \mi\right) + c_F^-F(l- \mi) \geq c_F
\end{align*}

Note that $l_5$ and $l_6$ are defined using minimization problems to account for both continuous and discrete probability distributions for the uncertain power supply. Therefore, the dispatch rule can be summarized as:

%\begin{tabular}{clllll}
%\hline
%case & $\pf$ & $\pii$ & $\pw$ & $\pv$ & if \\
%\hline
%1 & $0$ & $0$ & $\uw$ & $l-\uw$ & $0 < l < l_2$ \\
%2 & $0$ & $l-l_2$ & $\uw$ & $l_2-\uw$ & $l_2 \leq l \leq l_2+\mi$ \\
%3 & $0$ & $\mi$ & $\uw$ & $l-\mi-\uw$ & $l_2 + \mi \leq l < l_5$ \\
%4 & $l - l_7(l) - \mi$ & $\mi$ & $\uw$ & $l_7(l)-\uw$ & $l_5 \leq l \leq l_6$ \\
%5 & $\mf$ & $\mi$ & $\uw$ & $l-\mf-\mi-\uw$ & $l_6 \leq l$ \\
%\hline
%\end{tabular}
\vspace{4mm}
\begin{center}
{\footnotesize
\begin{tabular}{|l|lll|l|}
\hline
Rule \# & $\pw + \pv$ & $\pii$ & $\pf$ & applies if \\
\hline
1. & $l$ & $0$ & $0$ &  $0 \leq l \leq l_2$ \\
2. & $l_2$ & $l-l_2$ & $0$ &  $l_2 < l \leq l_2+\mi$ \\
3. & $l-\mi$ & $\mi$ & $0$ & $l_2 + \mi < l \leq l_5$ \\
4. & $l_7(l)$ & $\mi$ & $l - l_7(l) - \mi$ & $l_5 < l \leq l_6$ \\
5. & $l-\mf-\mi$ & $\mi$ & $\mf$  & $l_6 < l$ \\
\hline
\end{tabular}}
\end{center}

\vspace{5mm}
where
\begin{equation*}
l_7(l)  := \Fi{\dfrac{\cf - \pp{v-\cfu}\F{l-\mi-\mf} - \cfd\F{l-\mi}}{\cfu - \cfd}}
\end{equation*}

\hfill \Halmos \endproof

%The arbitrager's problem can be formulated as follows:
%
%\begin{subequations}
%\begin{align}
%& \off{\pv,\iv}{\lf\pv + \inw{\om}{}{\lb\iv}} \\
%& \acc{\pv + \iv = 0}
%\end{align}
%\end{subequations}
%

%For a given value of virtual generation $\pv$, the optimal forward dispatch is obtained through the following classic economic dispatch problem
%%
%\begin{subequations} \label{virtualb}
%\begin{align}
%& \of{\pw \geq 0, \pii \geq 0, \pf \geq 0}{\cf\pf + \ci\pii} \\
%& \acc{\pf + \pii + \pw + \pv = l: \lf} \\
%& \ac{\pw \leq \uw} \\
%& \ac{\pf \leq \mf} \\
%& \ac{\pii \leq \mi}
%\end{align}
%\end{subequations}
%%
%whose solution for $\ci < \cf$ is summarized in the table below.

%Likewise, for given values of $\pf,\pw,\pv$, the balancing re-dispatch actions are determined by solving:
%
%\begin{subequations}
%\begin{align}
%& \of{\pfu\geq 0,\pfd\geq 0, \iw, \sw \geq 0}{\cfu\pfu - \cfd\pfd + v\sw} \\
%& \acc{\pfu - \pfd + \iw + \sw - \pv = 0: \lb} \\
%& \ac{\pfd \leq \pf} \\
%& \ac{\pfu \leq \mf - \pf} \\
%& \ac{0 \leq \pw + \iw \leq \ww} \\
%& \ac{\sw \leq l} \\
%\end{align}
%\end{subequations}

\section{Proof of Corollary \ref{Corr:MO}} \label{Proof:MO}

\proof{Proof.}
The proof of point 1 follows directly from the comparison of the closed-form dispatch solutions provided in Theorems~\ref{Prop:ESCF} and~\ref{Prop:VB}: For $l_1 \geq l_2$ and $0 \leq l \leq min\left(\mi + l_1, l_5\right)$, dispatch rules 1--3 in these two theorems apply and provide exactly the same solution $(\pw, \pf, \pii)$.

Statement 2 is a consequence of the fact that if $\mf > \mw$, then $l_2 = F^{-1}\left(\frac{c_{I}}{c_{F}^{+}}\right)$ (see Corollary~\ref{Corr:CapacityAdequate}), and thus, $\nvs \geq \frac{c_I}{c_{F}^{+}}$ implies that $l_1 \geq l_2$. Furthermore, if $\mf > \mw$, it holds that $\mi + l_1 < l_5$, because substituting $\mi + l_1$ into~\eqref{L5} gives $\left(v-c_F^+\right)F\left(l_1-\mf\right) + c_F^+F(l_1) = \cfu \nvs \leq \cf$, since $\cfu \geq \cf$. From this point on, the proof of statement 2 proceeds as that of point 1.

Claim 3 is based on Corollary~\ref{Corr:FreeUW}. Indeed, the dispatch rule provided by this corollary is the same as that given by Theorem~\ref{Prop:VB} when $l \notin \left(l_5, \mw + \mi + \mf\right)$. Furthermore, note that $l_6 \leq \mi + \mf + \mw$, since substituting $l = \mi + \mf + \mw$ into~\eqref{L6} yields $\left(v-c_F^-\right)F\left(\mi + \mf + \mw-\mf- \mi\right) + c_F^-F(\mi + \mf + \mw- \mi) = \left(v-c_F^-\right)F\left(\mw\right) + c_F^-F( \mf + \mw) = v > \cf$.

Claim 4 stems from claim 3 by just noticing that $l_5 = \mw + \mi$ when $\cfu = \cf$ and $\mf > \mw$.

Statement 5 trivially follows from the fact that the dispatch solution provided by Theorem~\ref{Prop:VB} is one of the possible solutions that the stochastic dispatch rule admits when $\cf=\cfu=\cfd$ (see Corollary~\ref{Corr:FreeUW-DW}).
%The proof of this corollary is evident from the fact that, if $\frac{c_I}{c_{F}^{+}} \leq \frac{c_F-c_{F}^{-}}{c_F^{+}-c_{F}^{-}}$, the dispatch rule in Proposition~\ref{Prop:ESCF} can be obtained from that of Proposition~\ref{Prop:VB} by just replacing  $p_{W} + p_V$ in the latter with the optimal dispatch of the stochastic power capacity given by the efficient market.

Finally, the last claim relies on the fact that if $\cf=\cfd$, maximum cost-efficiency is always achieved by breaking the merit order (see Corollary~\ref{Corr:FreeDW}) except for $l \geq  l_4 + \mi$ when $l_3 > \mf$, and $l \geq \mf + \mi$ otherwise. Recall that the price-consistent conventional two-stage market always prompts dispatch solutions that respect the merit order.
\hfill \Halmos \endproof

\section{Proof of Theorem \ref{Prop:EJOR}} \label{Proof:EJOR}

\proof{Proof.} We divide the forward dispatch solutions that are feasible for model \eqref{EJOR} into two groups. The first group includes those solutions in which the inflexible generation is dispatched below its capacity ($\pii < \mi$). The solutions of the second group are characterized by the fact that the inflexible power generation technology is dispatched to its maximum capacity ($\pii = \mi$).

If $\pii < \mi$ and provided that $\ci < \cf$, the dispatch of the flexible generation $\pf$ has to be equal to 0 to comply with the optimality condition of the lower-level problem \eqref{EJOR_LL_OF}--\eqref{EJOR_LL_Final}. Or, in other words, the merit-order forward dispatch imposed through \eqref{EJOR_LL_OF}--\eqref{EJOR_LL_Final} implies that the more expensive flexible generation is only dispatched if the cheaper inflexible generation has reached its maximum capacity. For $\pf=0,\pii=l-\pw$ and using the expected balancing cost of Proposition \ref{Prop:Cost2Stage}, the bilevel problem \eqref{EJOR} reduces to the following single-level optimization problem:
\begin{subequations}
\begin{align}
& \of{\pw \geq 0}{\ci \pp{l-\pw} + v\inww{0}{\pw-\mf}{} + \cfu\inww{\pw-\mf}{\pw}{} } \\
& \acc{\pw \leq l}
\end{align}
\end{subequations}

Based on the KKT conditions of this convex optimization problem, the dispatch rule and system marginal expected cost in this case is given by:

\vspace{4mm}
\begin{center}
{\footnotesize
\begin{tabular}{|lll|l|l|}
\hline
$\pw$ & $\pii$ & $\pf$ & Marginal expected cost & applies if \\
\hline
$l$ & $0$ & $0$ & $(v-\cfu)\F{l-\mf}+\cfu\F{l}$ & $0 \leq l \leq l_2$ \\
$l_2$ & $l-l_2$ & $0$ & $\ci$ &  $l_2 \leq l \leq l_2 + \mi$ \\
\hline
\end{tabular}}\end{center}
\vspace{4mm}

%Since the forward dispatch of the stochastic generation is determined in the upper-level problem, and using Proposition XX for $\pf=0$, we can conclude that the marginal cost of the stochastic generation is $\pp{v-\cfu}\F{\pw-\mf}+\cfu\F{\pw}$. The stochastic generation would thus completely satisfy the load $l$ as long as $\pp{v-\cfu}\F{\pw-\mf}+\cfu\F{\pw}\leq \ci$, which implies that $l\leq l_2$. In that case, the marginal cost of the whole system is set by the stochastic generation. For values of load higher than $l_2$,

%Two possible strategies to be compared against each other. In the first one, the inflexible generation is dispatch below its capacity, i.e. $\pii < \mi$, and therefore the dispatch rule is:
%%
%\begin{equation}
%\begin{cases}
%\pw = l \qquad \pii = 0 \qquad \pf=0 & \ifq 0 \leq l \leq l_2 \\
%\pw = l_2 \qquad \pii = l-l_2 \qquad \pf=0 & \ifq l_2 \leq l \leq l_2 + \mi \\
%%\pw + \pf = \ph \qquad \pii = \mi & \ifq l_2+\mi \leq l \\
%\end{cases} \label{ejor1}
%\end{equation}
%
%The marginal cost of this dispatch rule is:
%
%\begin{equation}
%\begin{cases}
%(v-\cfu)\F{s-\mf}+\cfu\F{s} & \ifq 0 \leq l \leq l_2 \\
%\ci & \ifq l_2 \leq l \leq l_2 + \mi \\
%\end{cases}
%\end{equation}

Finally, the total expected cost corresponding to solutions with $\pii < \mi$, which we denote by $z^1(l)$, is computed as:

\begin{equation}
z^1(l) = \int_{0}^{\min(l,l_2)}A(s)ds + \ci \max(l-l_2,0)
\end{equation}

where $A(s) = (v-\cfu)\F{s-\mf}+\cfu\F{s}$.

On the other hand, the second group of feasible solutions are characterized by the fact that the inflexible generation technology is dispatched at full capacity, i.e., $\pii=\mi$. The dispatch rule of the flexible and stochastic generation can be thus derived from Proposition \ref{Prop:Cost2Stage} with $\pt=l-\mi$, that is,

\vspace{4mm}
\begin{center}
{\footnotesize
\begin{tabular}{|lll|l|l|}
\hline
$\pw$ & $\pii$ & $\pf$ & Marginal cost & applies if \\
\hline
$l-\mi$ & $\mi$ & $0$ & $(v-\cfu)\F{l-\mf}+\cfu\F{l}$ & $\mi \leq l \leq l_1+\mi$ \\
$l_1$ & $\mi$ & $l-l_1-\mi$ & $\cf + \pp{v-\cfu}\F{l-\mf}-\cfd\FF{l}$ &  $l_1+\mi \leq l \leq l_1+\mi+\mf$ \\
$l-\mi-\mf$ & $\mi$ & $\mf$ & $\pp{v-\cfd}\F{l-\mf}+\cfd\F{l}$ &  $l_1 + \mi + \mf \leq l$ \\
\hline
\end{tabular}}\end{center}
\vspace{4mm}

Consequently, the total expected cost in this second case, which we denote by $z^2(l)$, is given by:
\begin{equation}
z^2(l) = \ci\mi + \int_{0}^{\min(l-\mi,l_1)}A(s)ds + \int_{\min(l-\mi,l_1)}^{\min(l-\mi,l_1+\mf)}B(s)ds + \int_{\min(l-\mi,l_1+\mf)}^{l-\mi}C(s)ds
\end{equation}

\noindent where $B(s) = \cf + \pp{v-\cfu}\F{s-\mf}-\cfd\FF{s}$ and $C(s) = \pp{v-\cfd}\F{s-\mf}+\cfd\F{s}$.

Note that if $l \leq \mi$, the optimal forward dispatch must necessarily belong to the first group of feasible solutions and the total expected cost is, therefore, equal to $z^1(l)$. Likewise, if  $l \geq l_2+\mi$, the optimal dispatch requires $\pii=\mi$ and hence, the total expected cost is given by $z^2(l)$. For load levels $\mi \leq l \leq l_2+\mi$ we must, however, compare $z^1(l)$ against $z^2(l)$ to determine whether it is optimal to dispatch the inflexible generation at its maximum capacity or not. To conduct this comparison, we distinguish between two cases, namely, $l_2 \leq l_1$ and  $l_1 < l_2$.

If $l_2 \leq l_1$, and given that $\mi \leq l \leq l_2+\mi$, then $l - \mi \leq l_1$ and therefore:
\begin{equation}
z^2(l) = \ci\mi + \int_{0}^{l-\mi}A(s)ds
\end{equation}

Without any further assumption, we can rewrite $z^1(l)$ as:
\begin{equation}
z^1(l) = \ci \max(l-l_2,0) + \int_{0}^{l-\mi}A(s)ds + \int_{l-\mi}^{\min(l,l_2)}A(s)ds
\end{equation}

Therefore, for $\mi \leq l \leq l_2 + \mi$ and $l_2 \leq l_1$, we have:
\begin{equation}
z^2(l) - z^1(l) = \ci \pp{\mi-\max(l-l_2,0)} - \int_{l-\mi}^{\min(l,l_2)}A(s)ds
\end{equation}

Next we evaluate the function $z^2(l) - z^1(l)$ at the extremes of the interval $\mi \leq l \leq l_2+\mi$ and also compute its derivative with respect to $l$ and obtain the following results:

\begin{equation}
z^2(\mi) - z^1(\mi) = \begin{cases}
\ci l_2 - \int_{0}^{l_2}A(s)ds \geq 0 & \ifq l_2 \leq \mi \\
\ci \mi - \int_{0}^{\mi}A(s)ds \geq 0 & \ifq \mi < l_2 \\
\end{cases}
\end{equation}

\begin{equation}
z^2(l_2+\mi) - z^1(l_2+\mi) = 0
\end{equation}

\begin{equation}
\der{\left(z^2(l) - z^1(l)\right)}{l}= \begin{cases}
-\ci + A(l-\mi) \leq 0 & \ifq l_2 \leq l \\
- A(l) + A(l-\mi) \leq 0 & \ifq l < l_2 \\
\end{cases}
\end{equation}

\noindent where we have used that $A(s) \leq \ci$ for $s \leq l_2$ by definition, and that $A(s)$ is an increasing function. Hence, we can conclude that $z^1(l) \leq z^2(l)$ for $\mi \leq l \leq l_2 + \mi$ and provide the following dispatch rule for $l_2 \leq l_1$:

\vspace{4mm}
\begin{center}
{\footnotesize
\begin{tabular}{|lll|l|}
\hline
$\pw$ & $\pii$ & $\pf$ & applies if \\
\hline
$l$ & $0$ & $0$ & $0 \leq l \leq l_2$ \\
$l_2$ & $l-l_2$ & $0$ & $l_2 \leq l \leq l_2+\mi$ \\
$l-\mi$ & $\mi$ & $0$ & $l_2+\mi \leq l \leq l_1+\mi$ \\
$l_1$ & $\mi$ & $l-\mi-l_1$ & $l_1+\mi \leq l \leq l_1+\mi+\mf$ \\
$l-\mi-\mf$ & $\mi$ & $\mf$ & $l_1+\mi+\mf \leq l$ \\
\hline
\end{tabular}}\end{center}
\vspace{4mm}

Let us now analyze the case $l_1 < l_2$ for $\mi \leq l \leq l_2 + \mi$. First we rewrite $z^1(l)$ as:
\begin{align}
& z^1(l) = \ci \max(l-l_2,0)+\int_{0}^{\min(l-\mi,l_1)}A(s)ds + \int_{\min(l-\mi,l_1)}^{\min(l-\mi,l_1+\mf)}A(s)ds + \int_{\min(l-\mi,l_1+\mf)}^{l-\mi}A(s)ds + \nonumber \\
 & \qquad + \int_{l-\mi}^{\min(l,l_2)}A(s)ds
\end{align}

And therefore:
\begin{align}
& z^2(l)-z^1(l) = \ci\pp{\mi-\max(l-l_2,0)} - \int_{\min(l-\mi,l_1)}^{\min(l-\mi,l_1+\mf)}(A(s)-B(s))ds - \nonumber \\
& \qquad -  \int_{\min(l-\mi,l_1+\mf)}^{l-\mi}(A(s)-C(s))ds - \int_{l-\mi}^{\min(l,l_2)}A(s)ds
\end{align}

By evaluating $z^2(l)-z^1(l)$ at $\mi$, $l_2+\mi$ and $l_1+\mi$ we obtain:

\begin{equation}
z^2(\mi)-z^1(\mi) = \begin{cases}
\ci l_2 - \int_{0}^{l_2}A(s)ds \geq 0 & \ifq l_2 \leq \mi \\
\ci \mi - \int_{0}^{\mi}A(s)ds \geq 0 & \ifq \mi < l_2 \\
\end{cases}
\end{equation}

\begin{equation}
z^2(l_2+\mi)-z^1(l_2+\mi) = -\int_{l_1}^{\min(l_2,l_1+\mf)}(A(s)-B(s))ds - \int_{\min(l_2,l_1+\mf)}^{l_2}(A(s)-C(s))ds \leq 0
\end{equation}

\begin{equation}
z^2(l_1+\mi)-z^1(l_1+\mi) = \begin{cases}
\ci\mi - \int_{l_1}^{l_1+\mi}A(s)ds \geq 0 & \ifq l_1+\mi \leq l_2 \\
\ci\pp{l_2-l_1} - \int_{l_1}^{l_2}A(s)ds \geq 0 & \ifq l_2 < l_1+\mi \\
\end{cases}
\end{equation}

\noindent where we have used that $A(s) \geq B(s), \forall s\geq l_1$, that $A(s) \geq C(s), \forall s$, and that $A(s) \leq \ci, \forall s \leq l_2$. After checking that $\der{\left(z^2(l)-z^1(l)\right)}{l} \leq 0$ for $\mi \leq l \leq l_2+\mi$, we can conclude that \tred{there must exist at least one OR THERE MUST EXIST A UNIQUE??} $l_8$  such that $l_1 + \mi \leq l_8 \leq l_2 + \mi$ and $z^2(l_8)-z^1(l_8) = 0$. Therefore, the dispatch rule if $l_1 < l_2$ is:

\vspace{4mm}
\begin{center}
{\footnotesize
\begin{tabular}{|lll|l|}
\hline
$\pw$ & $\pii$ & $\pf$ & applies if \\
\hline
$l$ & $0$ & $0$ & $0 \leq l \leq \min(l_2,l_8)$ \\
$l_2$ & $l-l_2$ & $0$ & $l_2 \leq l \leq l_8$ \\
$l_1$ & $\mi$ & $l-\mi-l_1$ & $l_8 \leq l \leq l_1+\mi+\mf$ \\
$l-\mi-\mf$ & $\mi$ & $\mf$ & $\max(l_8,l_1+\mi+\mf) \leq l$ \\
\hline
\end{tabular}}\end{center}
\vspace{4mm}

\noindent where
\begin{align}
& l_8 := x: l_1 + \mi \leq x \leq l_2 + \mi \; \text{and} \; \ci\pp{\mi-\max(x-l_2,0)} - \int_{\min(x-\mi,l_1)}^{\min(x-\mi,l_1+\mf)}(A(s)-B(s))ds - \nonumber \\
& \qquad - \int_{\min(x-\mi,l_1+\mf)}^{x-\mi}(A(s)-C(s))ds - \int_{x-\mi}^{\min(x,l_2)}A(s)ds = 0
\end{align}

Consequently, the optimal dispatch rule prompted by the bilevel linear program~\eqref{EJOR} can be formulated as follows:

\vspace{4mm}
\begin{center}
{\footnotesize
\begin{tabular}{|c|lll|l|l|}
\hline
Rule \# & $\pw$ & $\pii$ & $\pf$ & \multicolumn{2}{c|}{applies if} \\
\hline
1. & $l$ & $0$ & $0$ && $0 \leq l \leq l_2$ \\
2. & $l_2$ & $l-l_2$ & $0$ && $l_2 \leq l \leq l_2+\mi$ \\
3. & $l-\mi$ & $\mi$ & $0$ & $l_2 \leq l_1$ & $l_2+\mi \leq l \leq l_1+\mi$ \\
4. & $l_1$ & $\mi$ & $l-\mi-l_1$ && $l_1+\mi \leq l \leq l_1+\mi+\mf$ \\
5. & $l-\mi-\mf$ & $\mi$ & $\mf$ && $l_1+\mi+\mf \leq l$ \\
\hline
6. &$l$ & $0$ & $0$ &\multirow{4}{2cm}{$l_1 < l_2$}& $0 \leq l \leq \min(l_2,l_8)$ \\
7. &$l_2$ & $l-l_2$ & $0$ && $l_2 \leq l \leq l_8$ \\
8. &$l_1$ & $\mi$ & $l-\mi-l_1$ && $l_8 \leq l \leq l_1+\mi+\mf$ \\
9. &$l-\mi-\mf$ & $\mi$ & $\mf$ & & $\max(l_8,l_1+\mi+\mf) \leq l$ \\
\hline
\end{tabular}}\end{center}
\vspace{4mm}

 \hfill \Halmos \endproof

\section{Proof of Corollary \ref{Corr:MO_EJOR}} \label{Proof:MO_EJOR}

\proof{Proof.}
Statement 1 is trivially inferred by comparing the tables provided in Theorems \ref{Prop:ESCF} and \ref{Prop:EJOR}: rules 1--5 in both theorems are identical. These rules apply for $l_1 \geq l_2$.

Statement 2 follows from the fact that any stochastic dispatch
solution that satisfies the merit order complies, by definition,
with the following two conditions simultaneously: i) it is an
optimal solution to the lower-level
problem~\eqref{EJOR_LL_OF}--\eqref{EJOR_LL_Final}, because it
respects the merit order \emph{and} ii) it minimizes the expected
system operating cost~\eqref{EJOR_OF}, because it is a solution
given by the stochastic dispatch rule.

Statement 3 simply highlights a particular case that is already
covered by statement 2, since any dispatch solution for which
$\pii = \mi$ preserves the merit order (recall that $\ci < \cf$).

Finally, statement 4 can be inferred from Corollary~\ref{Corr:FreeUW-DW} and by noticing that, if $\cf=\cfu=\cfd$, then $l_8 = l_2 + \mi$.
\hfill \Halmos \endproof

\end{APPENDICES}

%%
%\theendnotes

% Acknowledgments here
\ACKNOWLEDGMENT{The work of Juan M. Morales and Salvador Pineda is partly funded by DSF (Det Strategiske Forskningsr{\aa}d) through the CITIES (Juan M. Morales) and the 5's (Juan M. Morales and Salvador Pineda) projects, with identification numbers 1035-00027B and 12-132636, respectively.}

% References here (outcomment the appropriate case)

% CASE 1: BiBTeX used to constantly update the references
%   (while the paper is being written).
\bibliographystyle{ormsv080} % outcomment this and next line in Case 1
\bibliography{MeritOrder} % if more than one, comma separated

% CASE 2: BiBTeX used to generate mypaper.bbl (to be further fine tuned)
%\input{mypaper.bbl} % outcomment this line in Case 2

%If you don't use BiBTex, you can manually itemize references as shown below.

%%%%%%%%%%%%%%%%%
\end{document}
%%%%%%%%%%%%%%%%%